\magnification=1000
\hsize=11.7cm
\vsize=18.9cm
\lineskip2pt \lineskiplimit2pt
\nopagenumbers

\hoffset=-1truein
\voffset=-1truein

\advance\voffset by 4truecm
\advance\hoffset by 4.5truecm

\newif\ifentete

\headline{\ifentete\ifodd	\count0 
      \rlap{\head}\hfill\tenrm\llap{\the\count0}\relax
    \else
        \tenrm\rlap{\the\count0}\hfill\llap{\head} \relax
    \fi\else
\global\entetetrue\fi}

\def\entete#1{\entetefalse\gdef\head{#1}}
\entete{}

\input amssym.def
\input amssym.tex

\def\-{\hbox{-}}
\def\.{{\cdot}}
\def\O{{\cal O}}
\def\K{{\cal K}}
\def\F{{\cal F}}

\def\P{{\cal P}}

\def\G{{\cal G}}

\def\R{{\cal R}}

\def\H{{\cal H}}

\def\X{{\cal X}}

\def\ch{\frak c\frak h}

\def\Gr{\frak G\frak r}

\def\int{\frak i\frak n\frak t}

\def\qq{\quad{\rm and}\quad}

\def\mod{\frak m\frak o\frak d}

\def\too{\longrightarrow}
\def\aut{\frak a\frak u\frak t}

 3
 2
\font\large=cmr10  scaled \magstep 2
 2
 2
 2
\font\cds=cmr7

\centerline{\large On the reduction of Alperin's Conjecture}
\smallskip
\centerline{\large  to the quasi-simple groups}

\vskip 0.5cm

\centerline{\bf Lluis Puig }
\medskip
\centerline{\cds CNRS, Institut de Math\'ematiques de Jussieu}
\smallskip
\centerline{\cds 6 Av Bizet, 94340 Joinville-le-Pont, France}
\smallskip
\centerline{\cds puig@math.jussieu.fr}

\vskip 0.5cm
\noindent
{\bf £1. Introduction}
\bigskip
£1.1. In [10,~Chap.~16] we propose a refinement of Alperin's Conjecture whose possible proof can be reduced to check
that this refinement holds on the so-called {\it quasi-simple\/} groups. To carry out this checking obviously depends on admitting the 
{\it Classification of the Finite Simple Groups\/}, and our proof of the reduction itself uses the {\it solvability\/} of the
{\it outer automorphism group\/} of a finite simple group, a known fact whose actual proof depends on this classification [10,~16.11].

\medskip
£1.2.  Unfortunately, on the one hand  in [10,~Chap.~16] our proof also depends on checking, in the list of  {\it quasi-simple\/} groups, the technical condition [10,~16.22.1] --- which, as a matter of fact, is not always fulfilled: it is not difficult to exhibit a counter-example from Example~4.2 in [12]! -- and on the other hand, in the second half of [10,~Chap.~16] some arguments
 have been scratched. The purpose of this paper is to remove that troublesome condition and to repair the bad arguments there.{\footnote{\dag}{\cds In particular, this paper has to be considered as an ERRATUM of a partial contents of [10,~Chap.~16]
 going from [10,~16.20] to the end of the chapter.}} Eventually,  we find the better result stated below. 

\medskip
£1.3. Let us be more explicit. Let $p$ be a prime number, $k$ an algebraically closed 
 field of characteristic $p\,,$ $\O$ a complete discrete valuation ring of characteristic zero admitting $k$ as the {\it residue\/} field,
 $\hat G$ a $k^*\-$group of finite $k^*\-$quotient~$G$ [10,~1.23], $b$ a block of $\hat G$  [10,~1.25] and $\G_k (\hat G,b)$
 the {\it scalar extension\/} from $\Bbb Z$ to $\O$ of the {\it Grothendieck group\/} of the category of finite generated $k_*\hat Gb\-$modules [10,~14.3]. In [10,~Chap.~14], choosing a maximal Brauer $(b,\hat G)\-$pair $(P,e)\,,$ the existence of a suitable
 $k^*\-\Gr\-$valued functor $\widehat\aut_{(\F_{\!(b,\hat G)})^{^{\rm nc}}}$ over some full subcategory $(\F_{\!(b,\hat G)})^{^{\rm nc}}$ of the 
 {\it Frobenius $P\-$category\/}~$\F_{\!(b,\hat G)}$ [10,~3.2] allows us to consider an inverse limit of Grothen-dieck groups --- noted 
 $\G_k (\F_{\!(b,\hat G)},\widehat\aut_{(\F_{\!(b,\hat G)})^{^{\rm nc}}})$ and called the {\it Grothendieck group of\/} $\F_{\!(b,\hat G)}$ --- such that Alperin's Conjecture is actually equivalent to the existence of an $\O\-$module isomorphism [10,~I\hskip1pt 32 and Corollary~14.32]
 $$\G_k (\hat G,b)\cong \G_k (\F_{\!(b,\hat G)},\widehat\aut_{(\F_{\!(b,\hat G)})^{^{\rm nc}}})
 \eqno £1.3.1.$$

 \medskip
 £1.4. Denote by ${\rm Out}_{k^*}(\hat G)$ the group of {\it outer\/} $k^*\-$automorphisms of $\hat G$ and by 
${\rm Out}_{k^*}(\hat G)_b$ the stabilizer of $b$ in ${\rm Out}_{k^*}(\hat G)\,;$ it is clear that ${\rm Out}_{k^*}(\hat G)_b$
acts on $\G_k (\hat G,b)\,,$ and in [10,~16.3 and~16.4] we show that this group still acts on~$\G_k (\F_{\!(b,\hat G)},\widehat\aut_{(\F_{\!(b,\hat G)})^{^{\rm nc}}})\,.$ Our purpose is to show that the following statement
\smallskip
\noindent
(Q)\quad {\it For any $k^*\-$group with finite $k^*\-$quotient $G$ and any block $b$ of $\hat G\,,$ 
there is an $\O {\rm Out}_{k^*}(\hat G)_b\-$module isomorphism 
$$\G_k (\hat G,b)\cong \G_k (\F_{\!(b,\hat G)},\widehat\aut_{(\F_{\!(b,\hat G)})^{^{\rm nc}}})
\eqno £1.4.1\phantom{.}$$\/}
\smallskip
\noindent
can be proved by checking that it holds  in all the cases where $G$ contains a normal noncommutative simple subgroup 
$S$ such that $C_G (S) = \{1\}\,,$ $p$~divides~$\vert S\vert$ and~$G/S$ is a cyclic $p'\-$group.

\medskip
£1.5. To carry out this purpose, in [10,~Chap.~15] we develop reduction results relating both members of isomorphism~£1.4.1 with 
the Grothendieck groups coming from suitable  proper {\it normal sub-blocks\/}; recall that a  {\it normal sub-block\/} of~$(b,\hat G)$
is a pair $(c,\hat H)$ formed by a normal $k^*\-$subgroup $\hat H $ of~$\hat G$ and a block $c$ of $\hat H$ fulfilling $cb\not= 0\,.$ It is one of this reduction results --- namely in the case where $\hat G/\hat H$ is a $p'\-$group~[10,~Proposition~15.19] ---
 that we improve here, allowing us to remove condition [10,~16.22.1]. Then, following the same strategy as in [10,~Chap.~16],  we will show that in [10,~Chap.~16] all the statements including  this condition in their hypothesis can be replaced by stronger results and, at the end, we succeed in replacing [10,~Theorem~16.45] by the following more precise result.

\bigskip
\noindent
{\bf Theorem~£1.6.} {\it Assume that any block $(c,\hat H)$ having a normal sub-block~$(d,\hat S)$ of positive defect 
such that the $k^*\-$quotient $S$ of $\hat S$ is simple, $H/S$ is a cyclic $p'\-$group and $C_{H}(S) = \{1\}\,,$ fulfills the following two conditions
\smallskip
\noindent
{\rm £1.6.1}\quad ${\rm Out}(S)$ is solvable.
\smallskip
\noindent
{\rm £1.6.2}\quad There is an $\O {\rm Out}_{k^*}(\hat H)_c\-$module isomorphism
$$\G_k 
(\hat H,c)\cong\G_k(\F_{\!(c,\hat H)},\widehat \aut_{(\F_{\!(c,\hat H)})^{^{\rm nc}}})\,.$$
Then, for any block $(b,\hat G)$ there is an $\O {\rm Out}_{k^*} (\hat G)_b\-$module
isomorphism
$$\G_k (\hat G,b)\cong\G_k(\F_{\!(b,\hat G)},\widehat \aut_{(\F_{\!(b,\hat G)})^{^{\rm nc}}})
\eqno £1.6.3.$$\/}

\vfill
\eject

\bigskip
\bigskip
\noindent
{\bf £2. Notation and quoted results}
\bigskip

£2.1. We already have fixed $p\,,$ $k$ and $\O\,.$ We only consider $k\-$algebras $A$ of finite dimension
and denote by $J(A)$ the {\it radical\/} and by $A^*$ the group of invertible elements of~$A\,.$ Let $G$ be a finite group; 
a {\it $G\-$algebra\/}  is a $k\-$algebra $A$ endowed with a $G\-$action [4] and we denote by $A^G$ the subalgebra of 
$G\-$fixed elements. A  $G\-$algebra homomorphism from~$A$ to another   $G\-$algebra~$A'$ is a  {\it not necessarily unitary\/} $k\-$algebra homomorphism $f\,\colon A\to A'$ compatible
with the $G\-$actions; we say that $f$ is an {\it embedding\/} whenever  
$${\rm Ker}(f) = \{0\}\qq {\rm Im}(f) = f(1_A)A'f(1_A)
\eqno £2.1.1.$$

\medskip
£2.2. Recall that, for any subgroup $H$ of $G\,,$
a {\it point\/} $\alpha$ of $H$ on $A$ is an $(A^H)^*\-$conjugacy class
of primitive idempotents of $A^H$ and the pair $H_\alpha$ is a {\it
pointed group\/} on $A$ [5, 1.1]. For any $i\in \alpha\,,$ $iAi$ has an evident
structure of  $H\-$algebra and we denote by~$A_\alpha$ one of these mutually
$(A^H)^*\-$conjugate $H\-$algebras and by~$A(H_\alpha)$ the {\it simple
quotient\/} of $A^H$ determined by $\alpha\,.$ A second pointed group $K_\beta$ on $A$ is {\it contained\/} in~$H_\alpha$ if $K\i H$ and, 
for any $i\in\alpha\,,$ there is $j\in \beta$ such that~[5,~1.1]
$$ij = j = ji
\eqno £2.2.2;$$
then, it is clear that the $(A^K)^*\-$conjugation induces $K\-$algebra
embeddings
$$f_\beta^\alpha : A_\beta\too {\rm Res}^{H}_{K} (A_\alpha)
\eqno £2.2.3.$$

\medskip
£2.3. For any $p\-$subgroup $P$ of $G$ we consider the {\it Brauer quotient\/}
and the {\it Brauer homomorphism\/} [1, 1.2]
$${\rm Br}^A_P : A^P\too A (P) =  A^P\Big/\sum_Q A^P_Q
\eqno £2.3.1,$$
where $Q$ runs over the set of proper subgroups of $P\,,$ and call  {\it local\/} any point~$\gamma$ of $P$ 
on $A$ not contained in   ${\rm Ker(Br}^A_P)$ [5, 1.1]. Recall that
{\it a local pointed group $P_\gamma$ contained in $H_\alpha$ is maximal if and
only if ${\rm Br}_P(\alpha)\i A(P_\gamma)^{N_H (P_\gamma)}_P$\/} [5, Proposition~1.3] and then {\it the
$P\-$algebra~$A_\gamma$\/} --- called a {\it source algebra\/} of $A_\alpha$ --- {\it
is Morita equivalent to~$A_\alpha$\/} [9,~6.10]; moreover,
{\it the maximal local pointed groups~$P_\gamma$ contained in $H_\alpha$\/} ---
called the {\it defect pointed groups\/} of $H_\alpha$ --- {\it are  mutually
$H\-$conjugate\/} [5,~Theorem~1.2].

\medskip
£2.4. Let us say that $A$ is a {\it $p\-$permutation
$G\-$algebra\/} if a Sylow $p\-$subgroup of $G$ stabilizes a basis of $A$
[1, 1.1]. In this case, recall that if $P$ is a $p\-$subgroup of $G$ and $Q$ a normal subgroup of $P$
then the corresponding Brauer homomorphisms induce a $k\-$algebra isomorphism [1, Proposition~1.5]
$$\big(A(Q)\big)(P/Q)\cong A(P)
\eqno £2.4.1;$$
\eject
\noindent
moreover, choosing a point $\alpha$ of $G$ on~$A\,,$ we call 
{\it Brauer $(\alpha,G)\-$pair\/} any pair $(P,e_A)$ formed by a
$p\-$subgroup~$P$ of $G$ such that ${\rm Br}^A_P (\alpha)\not= \{0\}$ and by
 a primitive idempotent $e_A$ of the center $Z\big( A (P)\big)$ of $A(P)$ such
that 
$$e_A\. {\rm Br}^A_P (\alpha)\not= \{0\}
\eqno £2.4.2;$$
 note that any local pointed group $Q_\delta$ on $A$ {\it contained\/} in $G_\alpha$
determines a Brauer $(\alpha,G)\-$pair $(Q,f_A)$ fulfilling $f_A \.{\rm Br}^A_Q
(\delta)\not= \{0\}\,.$

\medskip
£2.5. It follows from [1,~Theorem~1.8] that {\it the
inclusion between the local pointed groups on $A$  induces an inclusion
between the Brauer $(\alpha,G)\-$pairs\/}; explicitly, if $(P,e_A)$ and
$(Q,f_A)$ are two Brauer $(\alpha,G)\-$pairs then we have
$$(Q,f_A)\i (P,e_A)
\eqno £2.5.1\phantom{.}$$
whenever there are local pointed groups $P_\gamma$ and $Q_\delta$ on $A$
fulfilling $$Q_\delta\i P_\gamma\i G_\alpha\quad ,\quad f_A \.{\rm Br}^A_Q (\delta)\not=
\{0\} \qq    e_A\. {\rm Br}^A_P (\gamma)\not= \{0\}
\eqno £2.5.2.$$
Actually, according to the same result, for any $p\-$subgroup $P$ of $G\,,$
any primitive idempotent $e_A$ of $Z\big(A(P)\big)$ fulfilling $e_A\. {\rm Br}^A_P (\alpha)\not= \{0\}$ 
and any subgroup~$Q$ of~$P\,,$ there is a unique
primitive idempotent $f_A$ of $Z\big(A(Q)\big)$ fulfilling 
$$e_A\. {\rm Br}^A_P (\alpha)\not=\{0\}\qq (Q,f_A)\i (P,e_A)
\eqno £2.5.3.$$
Once again, {\it the maximal Brauer $(\alpha,G)\-$pairs are pairwise
$G\-$conjugate\/} [1, Theorem~1.14].

\medskip
£2.6. For inductive purposes, we have to consider a  $k^*\-$group $\hat G$ of finite $k^*\-$quotient $G$ 
[10,~1.23] rather than a finite group; moreover, we are specially interested in the $G\-$algebras $A$ endowed with a $k^*\-$group
homomorphism $\rho\,\colon \hat G\to A^*$ inducing the action of $G$ on $A\,,$ called {\it $\hat G\-$interior algebras\/};
 in this case, for any pointed group $H_\alpha$ --- also noted $\hat H_\alpha$ --- on $A\,,$ $A_\alpha = iAi$ has a  structure of 
{\it $\hat H\-$interior algebra\/} mapping $\hat y\in\hat H$ on $\rho (\hat y)i = i\rho(\hat y)\,;$ moreover, setting 
$\hat x\.a\.\hat y = \rho(\hat x)a\rho(\hat y)$
for any $a\in A$ and any $\hat x,\hat y\in\hat  G\,,$ a  $\hat G\-$interior algebra homomorphism 
from~$A$ to another   $\hat G\-$interior algebra~$A'$ is a $G\-$algebra homomorphism $f\,\colon A\to A'$ fulfilling
$$f(\hat x\.a\.\hat y) = \hat x\.f(a)\.\hat y
\eqno £2.6.1.$$

\medskip
£2.7. In particular, if $H_\alpha$ and $K_\beta$ are two pointed groups on $A\,,$ we say that an
injective group homomorphism $\varphi\,\colon K\to H$ is an
{\it $A\-$fusion from $ K_\beta$ to~$H_\alpha$\/} whenever there is
a  $K\-$interior algebra~{\it embedding\/} 
$$f_{\varphi} : A_\beta\too {\rm Res}^{H}_{K} (A_\alpha) 
\eqno £2.7.1\phantom{.}$$
such that the inclusion $A_\beta\i A$ and the composition of $f_{\varphi}$
with the inclusion $A_\alpha\i A$ are $A^*\-$conjugate; we denote by $F_A
( K_\beta,H_\alpha)$ the set of $H\-$conjugacy classes of  $A\-$fusions from $ K_\beta$ to~$H_\alpha$
and, as usual, we write $F_A (H_\alpha)$ instead of $F_A
(H_\alpha,H_\alpha)\,.$ If $A_\alpha = iAi$ for $i\in \alpha\,,$  it follows
from [6, Corollary~2.13] that we have a group homomorphism
$$F_A (H_\alpha)\too N_{A_\alpha^{^*}}(H\.i)\big/H\.(A_\alpha^H)^*
\eqno £2.7.2\phantom{.}$$
and if $H$ is a $p\-$group then we consider the $k^*\-$group $\hat F_A (H_\alpha)$
defined by the {\it pull-back\/}
$$\matrix{F_A (H_\alpha)&\too& N_{A_\alpha^{^*}}(H\.i)/H\.(A_\alpha^H)^*\cr 
\uparrow&\phantom{\big\uparrow}&\uparrow\cr
\hat F_A (H_\alpha) &\too &N_{A_\alpha^{^*}}(H\.i)\big/H\.\big(i + J(A_\alpha^H)\big)\cr}
\eqno £2.7.3.$$

\medskip
£2.8. We also consider the {\it mixed\/} situation of an {\it $\hat H\-$interior $G\-$algebra~$B$\/} where $\hat H$ is a 
$k^*\-$subgroup of $\hat G$ and $B$ is a $G\-$algebra endowed with a {\it compatible\/} 
$\hat H\-$interior algebra structure, in such a way that the  $k_*\hat G\-$module $B\otimes_{k_*\hat H} k_*\hat G$ 
endowed with the product
$$(a\otimes\hat  x).(b\otimes\hat  y) = ab^{x^{-1}}\otimes\hat  x\hat y
\eqno £2.8.1,$$
for any $a,b\in B$ and any $\hat x,\hat y\in\hat  G\,,$ and with the group homomorphism mapping $\hat x\in\hat  G$ 
on $1_B\otimes\hat  x\,,$  becomes a $\hat G\-$interior algebra --- simply noted $B\otimes_{\hat H}\hat  G\,.$ 
For instance, for any $p\-$subgroup $P$ of $\hat G\,,$  $A(P)$ is a $C_{\hat G} (P)\-$interior $N_G (P)\-$algebra.

\medskip
£2.9. Obviously, the {\it group algebra\/} $k_*\hat  G$ is a $p\-$permutation $\hat G\-$interior algebra and, 
for any block $b$ of~$\hat G\,,$ the $\big((k_*\hat G)^G\big)^*\-$conjugacy class $\alpha = \{b\}$ is a {\it point\/} of $G$ on $k_*\hat  G\,.$ Moreover, for any $p\-$subgroup $P$ of $\hat G\,,$ the Brauer homomorphism 
 ${\rm Br}_P = {\rm Br}_P^{k_*\hat G}$ induces a $k\-$algebra isomorphism [8, 2.8.4]
$$k_*C_{\hat G} (P)\cong (k_*\hat  G)(P)
\eqno £2.9.1;$$
thus, up to identification throughout this isomorphism, in a Brauer $(\alpha,G)\-$ pair $(P,e)$  as defined above
--- called {\it Brauer $(b,\hat G)\-$pair\/} from now on --- $e$ is nothing but a block of $C_{\hat G} (P)$  
such that $e{\rm Br}_P (b)\not= 0\,.$ It is handy to consider the quotient
$$\bar C_{\hat G} (P) = C_{\hat G} (P)/Z(P)
\eqno £2.9.2\phantom{.}$$
and we denote by 
$${\rm \bar Br}_P : (k_*\hat G)^Q \too k_*\bar C_{\hat G} (P)
\eqno £2.9.3\phantom{.}$$
the corresponding homomorphism; recall that the image  $\bar e$ of $e$ in $k _*\bar C_{\hat G} (P)$ is a block of 
$\bar C_{\hat G} (P)$ and that  the {\it Brauer First Main Theorem\/} affirms that {\it $(P,e)$ is maximal if and only if the $k\-$algebra 
$k_*\bar C_{\hat G} (P)\bar e$  is simple and the inertial
quotient 
$$E_G (P,e) = N_{\hat G} (P,e)/P\.C_{\hat G} (P)
\eqno £2.9.4\phantom{.}$$
 is a $p'\-$group\/}~[9,~Theorem~10.14].

 \medskip
 £2.10. In this case, the {\it Frobenius $P\-$category\/}   $\F =\F_{\!(b,\hat G)}$ of $b$ [10,~3.2] is, up to equivalence, the category where the objects  are the Brauer $(b,\hat G)\-$pairs $(Q,f)$ and the morphisms are determined by the homomorphisms between the corresponding   $p\-$groups induced by the {\it inclusion\/} between Brauer $(b,\hat G)\-$pairs and by the $G\-$conjugation. 
 Then, we say~that $(Q,f)$  is {\it nilcentralized\/} if the block~$f$ of $C_{\hat G}(Q)$ is {\it nilpotent\/} 
 [10,~Proposition~7.2], we denote by 
 $\F^{^{\rm nc}}$ the full subcategory of $\F$ over the set of nilcentralized Brauer $(b,\hat G)\-$pairs, and
 consider the {\it proper category of $\F^{^{\rm nc}}\-$chains $\ch^*(\F^{^{\rm nc}})$\/} [10,~A2.8] and the 
 {\it automorphism functor\/} [10,~Proposition~A2.10]
 $$\aut_{\F^{^{\rm nc}}} : \ch^*(\F^{^{\rm nc}})\too  \Gr
 \eqno £2.10.1,$$
where $\Gr$ denotes the category of finite groups, mapping any $\ch^*(\F^{^{\rm nc}})\-$object $(\frak q,\Delta_n)$ --- $\frak q$ being a functor from the {\it ordered $n\-$simplex\/} $\Delta_n$ to $\F^{^{\rm nc}}$ --- to its $\ch^*(\F^{^{\rm nc}})\-$automorphism group --- the automorphism group of the functor
$\frak q\,,$ simply noted $\F(\frak q)\,.$

\medskip
£2.11. If $(Q,f)$ is a nilcentralized Brauer $(b,\hat G)\-$pair, $f$ determines a unique {\it local point\/} $\delta$ of $Q$ on 
$k_*\hat G$ since $k_*C_{\hat G}(Q)f$ has a unique isomorphism class of simple modules [7,~(1.9.1)]; now, the action 
of $N_G(Q,f)$ on the simple $k\-$algebra $(k_*\hat G)(Q_\delta)$ determines a $k^*\-$group $\hat N_G(Q,f)$ and it is clear that
the the corresponding $k^*\-$subgroup  $\hat C_G(Q)$ is canonically isomorphic to $C_{\hat G}(Q)\,,$ so that the
``difference'' $\hat N_G(Q,f) * N_{\hat G}(Q,f)^\circ$ [8,~5.9] admits a normal subgroup isomorphic to $C_G(Q)\,;$ then, up to identification, we define
$$\hat E_{G} (Q,f) = \big(\hat N_G(Q,f) * N_{\hat G}(Q,f)^\circ\big)\big/Q\.C_G(Q)
\eqno £2.11.1;$$
note that from [6, Theorem~3.1] and [8,~Proposition~6.12], suitable extended to $k^*\-$groups, we obtain {\it canonical\/} 
$k^*\-$group  isomorphisms (cf.~£2.8.3)
$$\hat E_G (Q,f)^\circ\cong \hat F_{k_*\hat G}(Q_\delta)\cong \hat F_{(k_*\hat G)_\delta}(Q_\delta)
\eqno £2.11.2.$$

\medskip
£2.12. In [10,~Theorem~11.32] we prove that the functor $\aut_{\F^{^{\rm nc}}}$ above can be lifted to a functor
$$\widehat{\aut}_{\F^{^{\rm nc}}} : \ch^*(\F^{^{\rm nc}})\too  k^*\-\Gr
 \eqno £2.12.1,$$
where $k^*\-\Gr$ denotes the category of $k^*\-$groups with finite $k^*\-$quotient, mapping any $\ch^*(\F^{^{\rm nc}})\-$object $(\frak q,\Delta_n)$ on the corresponding $k^*\-$subgroup $\hat\F (\frak q)$\break
\eject
\noindent
 of~$\hat E_{G} \big(\frak q (n)\big)\,.$
Then, denoting by $\frak g_k\,\colon k^*\-\Gr\to \O\-\mod$ the functor determined by the {\it Grothendieck groups\/} and 
the {\it restriction maps\/}, we define [10,~14.3.3]
$$\G_k(\F,\widehat{\aut}_{\F^{^{\rm nc}}}) = \lim_{\longleftarrow}\,(\frak g_k\circ \widehat{\aut}_{\F^{^{\rm nc}}})
\eqno £2.12.2$$
More precisely, we say that a Brauer $(b,\hat G)\-$pair is {\it $\F\-$selfcentralizing\/}  if the block~$\bar f$ of $\bar C_{\hat G}(Q)$
has {\it defect\/} zero [10,~Corollary~7.3]; denoting by $\F^{^{\rm sc}}$ the full subcategory of $\F$ over the set of selfcentralizing Brauer $(b,\hat G)\-$pairs and by $\widehat{\aut}_{\F^{^{\rm sc}}}$ the corresponding restriction, in [10,~Corollary~14.7] we prove that
$$\G_k(\F,\widehat{\aut}_{\F^{^{\rm nc}}})\cong \lim_{\longleftarrow}\,(\frak g_k\circ \widehat{\aut}_{\F^{^{\rm sc}}})
\eqno £2.12.3.$$

\medskip
£2.13. On the other hand, for any $p\-$subgroup $Q$ of~$\hat G$ and any $k^*\-$subgroup $\hat H$ of $N_{\hat G} (Q)$ 
containing $Q\.C_{\hat G}(Q)\,,$ we have
$${\rm Br}_Q\big((k_*\hat  G)^H\big) = (k_*\hat  G)(Q)^H
\eqno £2.13.1\phantom{.}$$
and therefore {\it any block $f$ of $C_{\hat G} (Q)$ determines a unique point~$\beta$ of $H$ on~$k_*\hat  G$
 such that $H_\beta$ contains $Q_\delta$ for a local point $\delta$ of $Q$ on $k_*G$
ful-filling\/}~[7,~Lemma~3.9]
$$f\. {\rm Br}_Q (\delta)\not= \{0\}
\eqno £2.13.2.$$
Recall that, if $R$ is a subgroup of $Q$ such that $C_{\hat G} (R)\i\hat  H$ then the blocks of $C_{\hat G} (R) 
= C_{\hat H}(R)$ determined by $(Q,f)$ from $\hat G$ and from $\hat H$ coincide
[1,~Theorem~1.8].

 \medskip
£2.14. Moreover, denote by $\gamma$ the local point of $P$  on $k_*\hat G$ determined by~$e$ and set 
$E_G (P,e) =E_G (P_\gamma)\,;$ since $E_G (P_\gamma)$ is a $p'\-$group, it follows from [9,~Lemma~14.10] that the short exact sequence
$$1\too P/Z(P) \too N_G (P_\gamma)/C_G (P)\too E_G (P_\gamma)\too 1
  \eqno £2.14.1\phantom{.}$$  
is split and that all the splittings are conjugate to each other; thus, any splitting determines an action of   $E_G (P_\gamma)$
 on $P$ and it is easily checked that the corresponding semidirect product
 $$\hat L_G(P_\gamma) = P\rtimes \hat E_G (P_\gamma)^\circ
 \eqno £2.14.2\phantom{.}$$
does not depend on our choice. Then,  it follows from  [9,~Theorem~12.8] that we have  a unique 
$\big((k_*\hat G)_\gamma^P\big)^*\-$conjugacy class of unitary $P\-$interior algebra homomorphisms
$$l_\gamma : k_*\hat L_G(P_\gamma)\too (k_*\hat G)_\gamma
\eqno £2.14.3\phantom{.}$$
which are also  $k(P\times P)\-$module {\it direct injections\/}.
\eject

\bigskip
\bigskip
\noindent
{\bf £3. Normal sub-blocks}
\bigskip

£3.1. Let~$\hat G$ be a $k^*\-$group with finite $k^*\-$quotient $G\,,$ $\hat H$ a normal $k^*\-$sub-group 
of $\hat G\,,$ $b$ a block of $\hat G$ and $c$ a block of $\hat H$ fulfilling  $cb\not= 0\,.$  Note that  we have 
$b{\rm Tr}_{\hat G_c}^{\hat G}(c) = b$ where $\hat G_c$ denotes the stabilizer of $c$ in $\hat G\,;$ thus, considering the $\hat G\-$stable
semisimple $k\-$subalgebra $\sum_{\hat x} k\.bc^{\hat x}$ of $k_*\hat  G\,,$ where $\hat x\in\hat  G$ runs over a set of representatives for $\hat G/\hat G_c\,,$ it follows from [11, Proposition~3.5] that $bc$ is a block of $\hat G_c$ and then from [11, Proposition~3.2] that we have
$$k_*\hat  Gb\cong {\rm Ind}_{\hat G_c}^{\hat G} (k_*\hat  G_c\, bc)
\eqno £3.1.1,$$
so that the source algebras of the block $b$ of $\hat G$ and of the block $bc$ of $\hat G_c$ are
isomorphic.

\medskip
£3.2. Thus, from now on we assume that $\hat G$ fixes~$c\,,$ so that we have $bc = b$  and, in particular, $\alpha = \{c\}$ is a point of $\hat G$ on $k_*\hat  H$ (cf.~£2.2). Let $(Q,f)$ be a maximal Brauer $(c,\hat H)\-$pair and denote by $N_{\hat G} (Q,f)$ the stabilizer of $(Q,f)$ in~$\hat G\,,$ setting 
$$C_{\hat G}(Q,f) = C_{\hat G}(Q)\cap N_{\hat G} (Q,f)
\eqno £3.2.1;$$
by the {\it Frattini argument\/}, we clearly get
$$\hat G = \hat H\.N_{\hat G} (Q,f)
\eqno £3.2.2;$$
as in~£2.11 above, $N_{G} (Q,f)$ acts on the simple $k\-$algebra $k_*\bar C_{\hat H}(Q)\bar f$ (cf.~£2.9), so that
 we get a $k^*\-$group $\hat N_{ G} (Q,f)$ and the ``difference'' $N_{\hat G} (Q,f) *\hat N_G (Q,f)^\circ$ contains
 a normal subgroup canonically isomorphic to $C_H(Q)\,;$ note that we have (cf.~£2.11.1)
 $$\hat E_H(Q,f)\i \big(N_{\hat G} (Q,f) *\hat N_G (Q,f)^\circ\big)\big/Q\.C_H (Q)
 \eqno £3.2.3.$$
 Moreover,  $C_{G} (Q,f)$ acts on the $k^*\-$group $\hat E_H (Q,f)$ acting
trivially on the $k^*\-$quotient $E_H (Q,f)\,,$ and therefore,  denoting by $S_{\hat G} (Q,f)$ the kernel of this action, the
quotient
$$Z = C_{\hat G} (Q,f)/S_{\hat G} (Q,f)
\eqno £3.2.4\phantom{.}$$
is an Abelian $p'\-$group.

\medskip
£3.3. More precisely, denoting by $\delta$ the local point of $Q$ on $k_*\hat H$ determined by~$f$ (cf.~£2.11) and choosing 
$j\in \delta\,,$ for any $\hat x\in N_{\hat G} (Q,f)$ there is $a_{\hat x}\in \big((k_*\hat H)^Q\big)^*$ such that 
$j^{\hat x} = j^{a_{\hat x}}$ and, in particular, the element $\hat x(a_{\hat x})^{-1}$ centralizes $j\,;$ hence,~choosing a set of
representatives $\X\i N_{\hat G} (Q,f)$ for the quotient $N_{\hat G} (Q,f)/N_{\hat H} (Q,f)\,,$ the element $\hat x(a_{\hat x})^{-1}$ normalizes the {\it source algebra\/} $B = j(k_*\hat H)j$ of $c$ for any $\hat x\in \X\,,$ and we easily get
$$ D = j(k_*\hat G)j = \bigoplus_{\hat x\in \X} \hat x(a_{\hat x})^{-1}\.B
\eqno £3.3.1.$$
\eject
\noindent

 It is clear that, for any $\hat x\in C_{\hat G} (Q,f)\,,$ the element  $\hat x(a_{\hat x})^{-1}$ induces a $Q\-$in-terior algebra automorphism of~$B$ and it follows from [8,~Proposition~14.9] that if $\hat x$ belongs to $S_{\hat G} (Q,f)$ then  the element 
$\hat x(a_{\hat x})^{-1}$ induces an {\it interior\/} automorphism of the  $Q\-$interior algebra~$B\,;$ thus, 
up to modifying our choice of $a_{\hat x}\,,$ we may assume that $\hat x(a_{\hat x})^{-1}$ centralizes~$B\,.$ From now on, we assume that $\hat x(a_{\hat x})^{-1}$ centralizes $B$ for any $\hat x\in  S_{\hat G} (Q,f)\,,$ so that $Z$ acts on $B\,,$ and that the elements of $\X$ belonging to  $C_{\hat G} (Q,f)\.N_{\hat H} (Q,f)$ are actually chosen in~$C_{\hat G} (Q,f)\,.$

\medskip
£3.4. On the other hand, denote by $\hat C_G(Q,f)$ the corresponding $k^*\-$sub-group of $\hat N_G(Q,f)$ and set
$$\hat C_{H}^{G} (Q,f) =  \big(C_{\hat G} (Q,f) *\hat C_G (Q,f)^\circ\big)\big/C_H (Q)
\eqno £3.4.1;$$
then, since $\bar f$ is a block of defect zero of $\bar C_{\hat H}(Q)\,,$ we have [11,~Theorem~3.7]
$$k_*\bar C_{\hat G} (Q,f)\bar f\cong k_*\bar C_{\hat H} (Q)\bar f \otimes_k k_*\hat C_{H}^{G} (Q,f)
\eqno £3.4.2;$$
more generally, we still have
$$k_*\bar C_{\hat G} (Q){\rm Tr}_{\bar C_{\hat G} (Q,f)}^{\bar C_{\hat G} (Q)} (\bar f) \cong
{\rm Ind}_{\bar C_{\hat G} (Q,f)}^{\bar C_{\hat G} (Q)} \big(k_*\bar C_{\hat H} (Q)\bar f \otimes_k k_*\hat C_{H}^{G} (Q,f)\big)
\eqno £3.4.3.$$

\medskip
£3.5. Note that, always from  [11,~Theorem~3.7], if $Q$ is a defect group of~$b$ then there is a block $\bar h$ of defect zero of 
$\hat C_{H}^{G} (Q,f)$ such that we have
$${\rm\bar Br}_Q(b) =  {\rm Tr}_{\bar N_{\hat G} (Q,\bar f\otimes \bar h)}^{\bar N_{\hat G} (Q)}(\bar f\otimes\bar h)
\eqno £3.5.1;$$
in this case, denoting by $S_{H}^{G} (Q,f)$ the image of  $S_{\hat G} (Q,f)$ in $C_{H}^{G} (Q,f)$ and by $\hat S_{H}^{G} (Q,f)$ 
the corresponding $k^*\-$subgroup  of~$\hat C_{ H}^{G} (Q,f)\,,$ it follows again from [11,~Theorem~3.7] that there is
a block $\bar \ell$ of defect zero  of $\hat S_{H}^{G} (Q,f)$ such that
$$\bar h\,  {\rm Tr}_{\hat C_{H}^{G} (Q,f)_{\bar \ell}}^{\hat C_{ H}^{G} (Q,f)}(\bar\ell) = \bar h\qq
k_*\hat C_{H}^{G} (Q,f)_{\bar \ell}\,\bar \ell = k_*\hat S_{H}^{G} (Q,f)\bar \ell\otimes_k k_*\hat Z_{\bar\ell}
\eqno £3.5.2$$
where  $\hat C_{H}^{ G} (Q,f)_{\bar \ell}$ is the stabilizer of~$\bar\ell$ in $\hat C_{H}^{ G} (Q,f)$ and $\hat Z_{\bar\ell}$ 
a suitable $k^*\-$group with the stabilizer $Z_{\bar\ell} = C_{H}^{ G} (Q,f)_{\bar \ell}/ S_{H}^{G} (Q,f)$ of~$\bar\ell$ in $Z$ (cf.~£3.3.1) as the $k^*\-$quotient.

\medskip
£3.6. Moreover,  it is clear that $E_H(Q,f)$ and $C_{H}^{ G} (Q,f)$ are normal subgroups of the quotient
$N_G (Q,f)/Q\.C_H (Q)$ and therefore their converse ima-ges $\hat E_H(Q,f)$ and~$\hat C_{H}^{ G} (Q,f)$ in the quotient 
$\big(N_{\hat G} (Q,f) *\hat N_G (Q,f)^\circ\big)\big/Q\.C_H (Q)$ (cf.~£3.2.3 and £3.4.1) still normalize each other; but, since we have
$$N_H(Q_\delta)\cap C_G (Q_\delta) = C_H (Q)\
\eqno £3.6.1,$$
 their commutator is contained in~$k^*\,;$ hence, $E_H(Q,f)$ also acts on  $\hat C_{H}^{ G} (Q,f)\,,$ acting trivially
  on $ C_{H}^{ G} (Q,f)$ and on $\hat  S_{H}^{G} (Q,f)\,.$ In particular, if $Q$ is a defect group of~$b$ then $E_H(Q,f)$ fixes~$\bar \ell$ and therefore
  it acts on the $k^*\-$group $\hat Z_{\bar\ell}\,,$ acting trivially on $Z_{\bar\ell}\,.$

\medskip  
£3.7.   But, the action of $C_G(Q,f)$  on $\hat E_G(Q,f)$ determines an injective group homomorphism (cf.~£3.2.4)
$$Z\too {\rm Hom}\big(E_H (Q,f),k^*\big)
\eqno £3.7.1.$$
Hence, since $Z$ is Abelian, if $Q$ is a defect group of~$b$ then  the action of  $E_H(Q,f)$ on  $\hat Z_{\bar\ell}$ induces a surjective group homomorphism
$$E_H (Q,f)\too {\rm Hom} (Z_{\bar\ell},k^*)
\eqno £3.7.2;$$
in this case, since $Z_{\bar\ell}$ is an Abelian $p'\-$group and we have 
$$Z(k_*\hat Z_{\bar\ell}) = k_*Z(\hat Z_{\bar\ell}) \cong k\check Z_{\bar\ell}
\eqno £3.7.3\phantom{.}$$
where $\check Z_{\bar\ell}$ denotes the image of $Z(\hat Z_{\bar\ell})$ in the $k^*\-$quotient of  $\hat Z_{\bar\ell}\,,$
the group $E_H(Q,f)$ acts transitively on the set of blocks of~$\hat Z_{\bar\ell}$ and, in particular, $\bar\ell$ is primitive in
$Z\big(k_*\hat C_{H}^{G} (Q,f)_{\bar \ell}\big)^{E_H(Q,f)}\,.$

\medskip
\noindent
{\bf Proposition~£3.8.} {\it With the notation above,  $b$ belongs to~$k_*\big(\hat H\.S_{\hat G}(Q,f)\big)$
and it is a block of $\hat H\.S_{\hat G}(Q,f)\,.$\/}
\medskip
\noindent
{\bf Proof:}  It follows from [10,~Proposition~15.10] that $b$  already belongs to $k_*\big(\hat H\.C_{\hat G}(Q,f)\big)$ and that it is a block of $\hat H\.C_{\hat G}(Q,f)\,;$ thus, with the notation above, we may assume that
$$\hat G = \hat H\.C_{\hat G}(Q,f)\qq C_{\hat H}(Q) = S_{\hat G}(Q,f)
\eqno £3.8.1;$$
in this case, since $\hat G/\hat H$ is a $p'\-$group, it follows from [10,~Proposition~15.9] that $Q$ is necessarily a defect group of $b\,.$

\smallskip
Consequently, since we have $S_{H}^{G}(Q,f) = \{1\}$ and $C_{\hat G}(Q,f) = C_{\hat G}(Q)$ (cf.~£3.8.1), it follows from~£3.7
above  that the unity element is primitive in the $k\-$algebra $Z\big(k_*\hat C_{H}^{G} (Q)\big)^{E_H(Q,f)}\,;$ but, we have (cf.~£3.4.2)
$$Z\big(k_*\bar C_{\hat G}(Q)\bar f \big)^{N_{\hat H}(Q,f)}\cong Z\big(k_*\hat C_{H}^{G} (Q)\big)^{E_H(Q,f)}
\eqno £3.8.2;$$
hence, the idempotent ${\rm\bar Br}_Q(c) = {\rm Tr}_{N_{\hat H}(Q,f)}^{N_{\hat H}(Q)}(\bar f)$ is also primitive in 
the $k\-$alge-bra $Z\big(k_*\bar C_{\hat G}(Q) \big)^{N_{\hat G}(Q)}\,,$ which forces ${\rm\bar Br}_Q(b) ={\rm\bar Br}_Q(c)\,.$
Since this applies to {\it any\/} block $b'$ of $\hat G$ such that $b'c = b'\,,$ we actually get $b = c\,.$

\medskip
£3.9. From now on, we assume that $\hat G/\hat H$ is a $p'\-$group; in particular, it follows from [10,~Proposition~15.9] that  $Q$ 
is necessarily a defect group of~$b\,;$ then, it follows from [10,~Lemma~15.16] that the local point $\delta$ of $Q$ on~$k_*\hat H$
in~£3.3 above splits into a set $\{(\delta,\varphi)\}_{\varphi\in \P(k_*\hat C_H^G(Q,f))}$ of local points of~$Q$ on~$k_*\hat G\,.$ Moreover,  the blocks $\bar h$ of $\hat C_H^G(Q,f)$ and $\bar\ell$ of $\hat S_H^G(Q,f)$ respectively determine  points 
$\varphi$ of $k_*\hat C_H^G(Q,f)$ and $\psi$ of $k_*\hat S_H^G(Q,f)\,,$ and it is quite clear that we have (cf.~£3.4.3)
$$(k_*\hat G)(Q_{(\delta,\varphi)})\cong {\rm Ind}_{\bar C_{\hat G} (Q,f)_\varphi}^{\bar C_{\hat G} (Q)} 
\Big(k_*\bar C_{\hat H} (Q)\bar f \otimes_k \big(k_*\hat C_{H}^{G} (Q,f)\big)(\varphi)\Big)
\eqno £3.9.1;$$
then, setting $Z_{\bar\ell} = Z_\psi\,,$ we also get a point $\theta$ of $k_*\hat Z_{\psi}$ such that (cf.~£2.5.2)
$$\big(k_*\hat C_{H}^{G} (Q,f)\big)(\varphi)\cong {\rm Ind}_{C_{H}^{G} (Q,f)_\psi}^{C_{H}^{G} (Q,f)} 
\Big(\big(k_*\hat S_H^G(Q,f)\big)(\psi) \otimes_k (k_*\hat Z_\psi)(\theta)\Big)
\eqno £3.9.2;$$
Denote by~$\check Z_\psi$ the image of  $Z(\hat Z_\psi)$ in $Z_\psi$  and consider the action of~$Z_\psi$ 
on~$B$ defined in~£3.3 above;  choosing an idempotent $i$ in the point $(\delta,\varphi)\,,$ our next result shows how to compute the {\it source algebra\/} $A = i(k_*\hat G)i$  of~$b$  from the $Q\-$interior algebra $B^{\check Z_\psi}\,.$

\bigskip
\noindent
{\bf Theorem~£3.10.} {\it  With the notation above, assume that $\hat G/\hat H$ is a $p'\-$group. Then the image of
$\hat E_H (Q_\delta)^{\check Z_\psi}$ in $E_H(Q_\delta)$ coincides with the intersection $E_H (Q_\delta)
\cap E_G (Q_{(\delta,\varphi)})\,,$ this equality can be lifted to a $k^*\-$group isomorphism
from $\hat E_H (Q_\delta)^{\check Z_\psi}$ to the converse image of this intersection
in $\hat E_G (Q_{(\delta,\varphi)})\,,$ and we have a $Q\-$interior algebra isomorphism
 $$A\cong B^{\check Z_\psi} \otimes_{\hat E_H (Q_\delta)^{\check Z_\psi}} 
 \hat E_G (Q_{(\delta,\varphi)})
 \eqno £3.10.1.$$.\/}

 \par \noindent
 {\bf Proof:} For any $\hat x,\hat y\in N_{\hat G}(Q,f)\,,$ it is clear that the ``difference'' between  
 $\hat x(a_{\hat x})^{-1} \hat y(a_{\hat y})^{-1}$ and $\hat x\hat y (a_{\hat x\hat y})^{-1}$ belongs to 
 $\big((k_*\hat H)^Q\big)^*\,,$ and therefore the union 
 $$\hat X = \bigcup_{\hat x\in N_{\hat G}(Q,f)} \hat x(a_{\hat x})^{-1}\.(B^Q)^*
  \eqno £3.10.2\phantom{.}$$
  is a $k^*\-$subgroup of $N_{D^*}(Q\.j)\,;$ thus, since we have
  $$(k_*\hat H)^Q\cap N_{\hat G}(Q,f) = C_{\hat H} (Q)
   \eqno £3.10.3,$$
    we obtain the exact sequence
 $$1\too (B^Q)^*\too \hat X\too N_{\hat G} (Q,f)/C_{\hat H}(Q)\too 1
 \eqno £3.10.4;$$
moreover, since $j$ is primitive in $B^Q\,,$ we still have the exact sequence
$$1\too Q\.\big(j + J(B^Q)\big)\too X\too N_{\hat G} (Q,f)/Q\.C_{\hat H}(Q)\too 1
\eqno £3.10.5\phantom{.}$$
where, as usual, $X$ denotes the $k^*\-$quotient of $\hat X\,.$
\eject

\smallskip
 Since the quotient $ \bar N = N_{\hat G} (Q,f)/Q\.C_{\hat H}(Q)$ is a $p'\-$group, this sequence is split and, actually, all the splittings are conjugate  [2,~Lemma~3.3 and~Proposition~3.5];  thus, denoting by $\skew3\hat{\bar N}$ the converse image in $\hat X$ 
 of a lifting   of~$\bar N$ to~$X\,,$ it is easily checked that  the $k^*\-$quotient  of $B\cap \skew3\hat{\bar N}$ is isomorphic to $E_H(Q,f)$ and therefore we may assume that (cf.~£2.14.3)
 $$B\cap \skew3\hat{\bar N} = l_\delta\big(\hat E_H(Q,f)\big)
 \eqno £3.10.6;$$
then, up to suitable identifications, isomorphism~£3.3.1  determines a $Q\-$interior algebra isomorphism
 $$D\cong B\otimes_{\hat E_H(Q,f)} \skew3\hat{\bar N}
  \eqno £3.10.7.$$

  \smallskip
Moreover, identifying $C_H^G (Q,f)$ with its image in $\bar N\,,$ it is clear that $C_H^G (Q,f)$ centralizes $Q\.j\,;$
further, the action on $B$ of an element of $K_H^G (Q,f)$ coincides with the action of some element in $ j + J(B^Q)$ and thus
$K_H^G (Q,f)$ acts trivially on $B\,.$ On the other hand, since the group of fixed points $\bar N^Q$ of $Q$ on  
$\bar N$ coincides with $C_H^G (Q,f)$ and since from isomorphism~£3.10.7 we clearly get
$$D(Q) \cong B(Q)\otimes_{\hat E_H(Q,f)} \skew3\hat{\bar N}^Q
\eqno £3.10.8,$$
it follows from isomorphism~£3.4.2 that we also get a $k^*\-$isomorphism
$$\skew3\hat{\bar N}^Q\cong \hat C_H^G (Q,f)
\eqno £3.10.9.$$

\smallskip 
Firstly  assume that $\hat G = \hat H\.S_{\hat G}(Q,f)\,;$ in this case, we have $Z = \{1\}$ and, according to~£3.6,  isomorphism~£3.10.8 above becomes
$$D\cong B\otimes_k k_*\hat K_H^G(Q,f)
\eqno £3.10.10;$$
hence, we may assume that $i = j\otimes \ell$ for some primitive idempotent $\ell$ in the $k\-$algebra 
$k_*\hat K_H^G(Q,f)\,;$ then, $i$ centralizes $B$ and the multiplication by $i$ determines a $Q\-$interior algebra 
isomorphism $B \cong A\,;$ in particular, we get (cf.~£2.11.2)
$$\hat E_H(Q_\delta) \cong \hat F_B(Q_\delta)^\circ\cong \hat F_A(Q_{(\delta,\varphi)})^\circ
\cong \hat E_G(Q_{(\delta,\varphi)})
\eqno £3.10.11.$$

\smallskip
Consequently, in order to prove the theorem  we may assume that $\hat H$ contains $S_{\hat G}(Q,f)$ and,
in this case, firstly assume that $\hat G = \hat H\.C_{\hat G}(Q,f)\,;$ then, we have $K_H^G(Q,f) = \{1\}\,,$ $\bar \ell = 1$
and $\psi = \{1\}\,,$ and isomorphism~£3.10.7 above becomes (cf.~£3.6)
$$D\cong B\otimes_{k^*} \hat Z
\eqno £3.10.12;$$
in particular, $D^Q$ contains the $k\-$algebra $k_*\hat Z\,;$ moreover, since $\hat E_H(Q,f)$ and $\hat Z$ normalize each other
(cf.~£3.6), $\hat E_H(Q,f)$ normalizes the $k\-$subalgebra $k_*\hat Z$ of $D$ and, according to~£3.7 above,
it acts transitively on the set of primitive\break
\eject
\noindent
 idempotents of~$Z(k_*\hat Z)\,;$ but, it is clear that we have
$Z(k_*\hat Z) = k_*Z(\hat Z)$ and that its primitive idempotents have the form
$$e_\theta = {1\over \vert\check Z\vert}\.\sum_{\hat w}\theta(\hat w)\.\hat w^{-1}
\eqno £3.10.13\phantom{.}$$
where $\hat w\in Z(\hat Z)$ runs over a set of representatives for $\check Z$ and $\theta\,\colon Z(\hat Z)\to k^*$
is a $k^*\-$group homomorphism; hence, choosing such a  $k^*\-$group homomorphism~$\theta\,,$ it follows from [11,~Proposition~3.2] that we have
$$D\cong {\rm Ind}_{\hat E_H(Q,f)_\theta}^{\hat E_H(Q,f)} (C\otimes_{k} k_*\hat Ze_\theta)
\eqno £3.10.14\phantom{.}$$
where $\hat E_H(Q,f)_\theta$ denotes the stabilizer of $\theta$ in $\hat E_H(Q,f)$ and  $C$ the centralizer of the simple $k\-$algebra $k_*\hat Ze_\theta$ in $e_\theta D e_\theta\,.$

\smallskip
On the one hand, since the action of $\hat E_H(Q,f)$ on $\hat Z$ determines a homorphism from $E_H(Q,f)$ to the group ${\rm Hom}(Z,k^*)$
which is Abelian, $\hat E_H(Q,f)_\theta$ is normal in $\hat E_H(Q,f)$ and therefore $\hat E_H(Q,f)_\theta$ coincides with
$\hat E_H(Q,f)^{\check Z}\,.$ On the other hand, since $Z$ is a $p'\-$group, an elementary computation shows that
$$e_\theta (B\otimes_{k^*} \hat Z) e_\theta = B^{\check Z}\otimes_k k_*\hat Ze_\theta
\eqno £3.10.15\phantom{.}$$
and therefore we get $C= B^{\check Z}\,.$ In particular, since the unity element $j$ is primitive in $(B^{\check Z})^Q\,,$ 
 up to suitable identifications, in~isomorphism~£3.10.14 above we may assume that $\varphi = \theta$ and 
$i = 1\otimes (j\otimes e_\theta)\otimes 1\,,$ so that we obtain a $Q\-$interior algebra isomorphism $A\cong B^{\check Z}\,;$ moreover, once again because of~$Z$ is a $p'\-$group, we get (cf.~£2.11.2)
$$\hat E_H(Q_\delta)^{\check Z} \cong \big(\hat F_B(Q_\delta)^\circ\big)^{\check Z}\cong \hat F_A(Q_{(\delta,\varphi)})^\circ
\cong \hat E_G(Q_{(\delta,\varphi)})
\eqno £3.10.16.$$

\smallskip
Finally, in order to prove the theorem  we may assume that $\hat H$ contains~$C_{\hat G}(Q,f)\,;$ then, we have 
$K_H^G(Q,f) = \{1\} = C_H^G(Q,f)\,,$ $Z =\{1\}\,,$ $\bar \ell = 1 = \bar h$ and $\psi = \{1\} =\varphi\,,$ and in particular we get 
a group isomorphism
$$\bar N = N_{\hat G}(Q,f)/Q\.C_{\hat H}(Q,f) \cong E_G(Q_{(\delta,\{1\})})\eqno £3.10.17.$$
In this case we claim that $i = j\,;$ indeed, it  is clear that the multiplication by $B$ on the left and the action
 of~$Q$ by conjugation endows $D$ with a $B\rtimes Q\-$module structure and, since the idempotent $j$ is primitive in $B^Q\,,$ equality~£3.3.1 provides  a direct sum decomposition of $D$ on $B\rtimes Q\-$modules. More explicitly, note that $B$ is an indecomposable $B\rtimes Q\-$module since we have
 ${\rm End}_{B\rtimes Q}(B) = B^Q\,;$ but,
for any $\hat x\in \X\,,$ the inversible element $\hat x(a_{\hat x})^{-1}j$ of $D$ together with the action of $\hat x$ on~$Q$ determine an automorphism $g_{\hat x}$ of~$B\rtimes Q\,;$ thus,   equality~£3.3.1 provides the following direct sum decomposition on indecomposable $B\rtimes Q\-$modules
 $$D\cong \bigoplus_{\hat x\in  \X} {\rm Res}_{g_{\hat x}}(B)
 \eqno £3.10.18.$$

 \smallskip
 Moreover, we claim that the $B\rtimes Q\-$modules ${\rm Res}_{g_{\hat x}}(B)$ and ${\rm Res}_{g _{\hat x'}}(B)$ for 
 $\hat x,\hat x'\in \X$  are isomorphic if and only if $\hat x = \hat x'\,;$ indeed, a $B\rtimes Q\-$module isomorphism
 $${\rm Res}_{g_{\hat x}}(B)\cong {\rm Res}_{g_{\hat x'}}(B)
 \eqno £3.10.19\phantom{.}$$
 is necessarily determined by the multiplication on the right by an inversible element $a$ of $B$ fulfilling
 $u^{\hat x}\.a = a\.u^{\hat x'}$ or, equivalently, $(u\.j)^a = u^{\hat x^{-1}\hat x'}\.j$ for any~$u\in Q\,,$ which amounts to saying that the automorphism of $Q$ determined by  $\hat x^{-1}\hat x'\in N_{\hat G}(Q,f)$ is a {\it $B\-$fusion\/} (cf.~£2.7) from 
 $Q_\delta$  to $Q_\delta$ [6,~Proposition~2.12]; but, it follows from [6, Proposition~2.14 and Theorem~3.1] that we have
$$ F_B(Q_\delta) = E_H (Q,f)
\eqno £3.10.20\,;$$
then, isomorphism~£3.10.17 implies that $\hat x^{-1}\hat x'$ belongs to $N_{\hat H}(Q,f)\,,$ so that we still have  
$\hat x = \hat x'\,.$

\smallskip
On the other hand, it is clear that $Di$ is a direct summand of $D$ as $B\rtimes Q\-$modules and therefore there is 
$\hat x\in \X$ such that ${\rm Res}_{g_{\hat x}}(B)$ is a direct summand of the $B\rtimes Q\-$module $Di\,;$ but, it follows from
[6, Proposition~2.14] that we have
$$F_D (Q_{(\delta,\{1\})}) = F_A (Q_{(\delta,\{1\})}) = E_G(Q_{(\delta,\{1\})})
\eqno £3.10.21\phantom{.}$$
and therefore, once again applying [6,~Proposition~2.12], for any element $\hat y$ in $N_{\hat G}(Q_{(\delta,\{1\})})$
there is an inversible element $d_{\hat y}$ in $D$ fulfilling  
$$(u\.i)^{d_{\hat y}} = u^y\.i
\eqno £3.10.22\phantom{.}$$
for any $u\in Q\,;$ then, for any $\hat x'\in \X\,,$ it is clear that $Di = Did_{\hat x^{-1}\hat x'}$ has a direct summand
isomorphic to  ${\rm Res}_{g_{\hat x'}}(B)\,,$ which forces the equality of the dimensions of $Di$ and $D\,,$ proving our claim.

\smallskip
Consequently, from isomorphism~£3.10.17 the  $Q\-$interior algebra isomorphism~£3.10.7 becomes$$A\cong B
\otimes_{\hat E_H(Q_\delta)} \widehat{E_G(Q_{(\delta,\{1\})})}
\eqno £3.10.23,$$
for a suitable $k^*\-$group $\widehat{E_G(Q_{(\delta,\{1\})})}$ with $k^*\-$quotient $E_G(Q_{(\delta,\{1\})})\,,$
and then it easily follows from~£2.14.1 that we have
$$\hat E_H(Q_\delta)\i  \widehat{E_G(Q_{(\delta,\{1\})})}\cong \hat E_G(Q_{(\delta,\{1\})})
\eqno £3.10.24.$$
We are done.

\medskip
£3.11. As a matter of fact, this theorem implies [10, Corollary~15.20] {\it without assuming condition\/} [10, 15.17.1],
 as we show in the next result.
 \eject

\bigskip
\noindent
{\bf Corollary~£3.12.} {\it  With the notation above, assume that  $\hat G/\hat H$ is a $p'\-$group.
Let $\hat G'$ be a $k^*\-$subgroup of $\hat G$ containing $\hat H\,,$ $b'$ a block of~$\hat G'$~such that
 ${\rm Br}_Q(b')\not= 0\,,$ and $\varphi'$ a point of the   $k\-$algebra $k_*  C_{\hat H}^{\hat G'}(Q,f)$ 
 such that  $Q_{(\delta,\varphi')}$ is a defect pointed group of~$b'\,.$ If we have $E_{G'} (Q_{(\delta,\varphi')}) = 
 E_{G} (Q_{(\delta,\varphi)})$ and this equality can be lifted to a $k^*\-$group isomorphism
 $\hat E_{G'} (Q_{(\delta,\varphi')}) \cong \hat E_{G} (Q_{(\delta,\varphi)})$ then the
Frobenius $Q\-$category $\F' = \F_{\!(b',\hat G')}$ coincides with $\F\,,$ we have a
natural isomorphism $\widehat\aut_{\F'^{^{\rm nc}}}\cong \widehat\aut_{\F^{^{\rm nc}}}$
inducing an $\O\-$module isomorphism
$$\G_k (\F,\widehat\aut_{\F^{^{\rm nc}}})\cong \G_k (\F',\widehat\aut_{\F'^{^{\rm nc}}})
\eqno £3.12.1,$$
and the restrictions to the respective source algebras induce an $\O\-$module
isomorphism 
$$\G_k (\hat G,b) \cong \G_k (\hat G',b')
\eqno £3.12.2.$$\/}

\par
\noindent 
{\bf Proof:} Since we assume that $E_{G'} (Q_{(\delta,\varphi')}) = E_{G} (Q_{(\delta,\varphi)})\,,$ we have
$$E_H (Q_\delta)\cap E_{G'} (Q_{(\delta,\varphi')}) = E_H (Q_\delta)\cap E_{G} (Q_{(\delta,\varphi)})
\eqno £3.12.3\phantom{.}$$
and therefore it follows from Theorem~£3.10 that we still have a $k^*\-$group isomorphism
$$\hat E_H (Q_\delta)^{\check Z'_{\psi'}}\cong \hat E_H (Q_\delta)^{\check Z_\psi}
\eqno £3.12.4\phantom{.}$$
which actually forces $\check Z'_{\psi'}\cong \check Z_{\psi}$ (cf.~£3.6 and~£3.7) and $B^{\check Z'_{\psi'}}\cong 
B^{\check Z_{\psi}}\,.$ Thus , denoting by $A'$ a {\it source algebra\/} of the block $b'\,,$ always from Theorem~£3.10
we obtain $Q\-$interior algebra isomorphisms
$$\matrix{A'&\cong &B^{\check Z'_{\psi'}}\otimes_{\hat E_H (Q_\delta)^{\check Z'_{\psi'}}} \hat E_{G'} (Q_{(\delta,\varphi')})\cr
&&\wr\Vert\cr
&&B^{\check Z'_{\psi'}}\otimes_{\hat E_H (Q_\delta)^{\check Z_{\psi}}} \hat E_{G} (Q_{(\delta,\varphi)})& \cong &A\cr}
\eqno £3.12.5.$$
But, it follows from [6,~Theorem~3.1] and from [8,~Proposition~6.21] that  $\F$ and $\F'\,,$ $\widehat\aut_{\F^{^{\rm nc}}}$ and $\widehat\aut_{\F'^{^{\rm nc}}}\,,$
$\G_k (\F,\widehat\aut_{\F^{^{\rm nc}}})$ and $\G_k (\F',
\widehat\aut_{\F'^{^{\rm nc}}})\,,$ and $\G_k (\hat G,b)$ and $\G_k (\hat G',b')$ are
completely determined from the respective source algebras $A$ of $b$ and $A' $ of $b'\,.$ Thus, the isomorphism $A\cong A'$
forces the equality $\F = \F'$ and all the isomorphisms. We are done.

 \bigskip
\bigskip
\noindent
{\bf £4. Reduction of the question (Q)}
\bigskip
£4.1. From now on, we prove Theorem~£1.6 by revising all the contents of~[10,~Chap.~16]. The point is that there all the reduction
arguments depend on condition [10,~16.22.1] only thoughout condition  [10, 15.17.1] in  [10, Corollary~15.20]; since
this condition has been removed in~Corollary~£3.12 above, obtaining the same conclusion, it is possible to remove 
 condition [10,~16.22.1]\break
\eject
\noindent in all the statements of [10,~Chap.~16], proving Theorem~£1.6. We revise step by step,
 avoiding as far as possible to repeat proofs in the first part of the proof; from £4.11 on, we have to replace the corresponding
 part in [10,~Chap.~16] by new arguments.

 \medskip
 £4.2. Let $\hat G$ be a $k^*\-$group with finite $k^*\-$quotient $G$ and $b$ a block of~$\hat G\,;$
 from [10,~Proposition~16.6] we may assume that, for any nontrivial charac-teristic $k^*\-$subgroup  $\hat N$ of  $\hat G\,,$ any block~$c$ of~$\hat N$ such that $bc\not= 0$ is $\hat G\-$stable; then, from  [10,~Proposition~16.7] we may assume that,
for any nontrivial characteristic $k^*\-$subgroup  $\hat N$ of  $\hat G\,,$ any block~$c$ such that $bc = b$ has a nontrivial defect group.

\medskip
£4.3. From now on, we assume that, for any nontrivial characteristic $k^*\-$subgroup  $\hat N$ of  $\hat G\,,$ any block~$c$ of~$\hat N$ such that $bc\not= 0$ is $\hat G\-$stable and  has a nontrivial defect group, which forces $\Bbb O_{p'}(\hat G) = k^*\,.$
Then, from [10,~Proposition~£16.8] we may assume  that the quotient $\hat G/C_{\hat G}\Big(Z\big(\Bbb O_{p}(\hat G)\big)\Big)$
is a cyclic
 $p'\-$group and, moreover, from [10,~Proposition~£16.9] we may assume that we actually have $\Bbb O_{p}(\hat G) = \{1\}\,.$
Consequently, from now on we also assume that
$$\Bbb O_{p'}(\hat G) = k^*\qq \Bbb O_{p}(\hat G) = \{1\}
\eqno £4.3.1.$$

\medskip
£4.4. Then, it is well-known that the product $H$ of all the minimal non-trivial normal
subgroups of $G$ is a {\it characteristic\/} subgroup
of~$G$ isomorphic to a direct product 
$$H\cong \prod_{i\in I} H_i
\eqno £4.4.1\phantom{.}$$
of a finite family of noncommutative simple groups $H_i$ of order divisible by~$p$ [3,~Theorem~1.5].
 Denoting by $\hat H$ and by $\hat H_i$ the respective converse
images of $H$ and $H_i$ in $\hat G\,,$ it is quite clear that
$$\hat H = \widehat{\prod}_{i\in I} \hat H_i\quad{\rm and}\quad
C_{\hat G} (\hat H) = k^*
\eqno £4.4.2\phantom{.}$$
where $\,\widehat{\prod}_{i\in I}\hat H_i$ denotes the obvious {\it central product\/} of the family of  $k^*\-$groups $\hat H_i$
over~$k^*\,;$ moreover, since this decomposition is unique, the action of ${\rm Aut}_{k^*}(\hat G)$ on
$\hat H$ induces an ${\rm Aut}_{k^*}(\hat G)\-$action on~$I$ and, denoting by $\hat W$
 the kernel of the action of $\hat G$ on~$I\,,$ we have  $\hat H\i\hat W$ and get an
injective group homomorphism
$$\hat W/\hat H\too \prod_{i\in I}{\rm Out}_{k^*} (\hat H_i)
\eqno £4.4.3;$$
thus, admitting the announced {\it Classification of the Finite Simple
Groups\/}, the quotient $\hat W/\hat H$ is {\it solvable\/}.

\medskip
£4.5. Let $c$ be the block of $\hat H$ such that $cb = b$ and $(P,e)$ a maximal Brauer $(b,\hat G)\-$pair; 
setting $Q = P\cap \hat H\,,$ it follows from [10,~Proposition~15.9] that $Q$ is a defect group of $c$ and
that there is a block $f$ of $C_{\hat H}(Q)$ such that we have $e{\rm Br}_P (f)\not= 0$ and that
 $(Q,f)$ is a maximal Brauer $(c,\hat H)\-$pair; then, we consider the Frobenius $P\-$
and $Q\-$categories [10,~3.2]
$$\F = \F_{\!(b,\hat G)}\qq \H = \F_{\!(c,\hat H)}
\eqno £4.5.1.$$
Since clearly $c = \otimes_{i\in I} c_i$ where $c_i$ is a block of $\hat H_i\,,$ we have  
$Q =\prod_{i\in I} Q_i$ where~$Q_i$ is a defect group of $c_i\,,$ and $f = \otimes_{i\in I}
f_i$ where $f_i$ is a	block of~$C_{\hat H_i}(Q_i)$ and $(Q_i,f_i)$ is a maximal Brauer
$(c_i,\hat H_i)\-$pair.

\medskip
£4.6. Moreover, since we are assuming that any block involved in $b$ of~any
nontrivial characteristic $k^*\-$subgroup of~$\hat G$ has positive defect, for any~$i\in I$
the defect group $Q_i$ is {\it nontrivial\/}; thus, since any
{\it $\H\-$selfcentralizing\/} subgroup $T$ of $Q$ [10,~4.8] contains $Z(Q) = \prod_{i\in I} Z(Q_i)\,,$ 
$C_{\hat G}(T)$ centralizes
 $Z(Q_i)\not=\{1\}$ for any $i\in I$ and therefore we get
$$C_{\hat G}(T)\i \hat W
\eqno £4.6.1.$$
In particular, $W$ contains $\hat K = \hat H\.C_{\hat G}(Q,f)\,,$ which is actually a normal subgroup of $\hat G$ by the {\it Frattini argument\/}, and therefore the quotient $\hat K/\hat H$ is solvable (cf.~£16.11). Then, from [10,~Proposition~16.15] we may assume that this quotient has
{\it $p\-$solvable length\/} $1$ and that $\hat G/\hat K$ is a cyclic $p'\-$group; going further,
from  [10,~Proposition~16.19] we actually may assume that $\hat G/\hat H$ is a $p'\-$group and that $\hat G/\hat K$ is cyclic.

\medskip
£4.7. Consequently, from now on we assume that $\hat G/\hat H$ is a $p'\-$group and that $\hat G/\hat K$ is cyclic. 
At this point, since Corollary~£3.12 holds, we can remove condition~[10,~16.22], and then Proposition~16.23 in~[10]
becomes.

\bigskip
\noindent
{\bf Proposition~£4.8}\phantom{.} {\it With the notation above, assume that $\hat G/\hat H$ is a $p'\-$group
and that $C = \hat G/\hat K$ is cyclic. Denote by
$\delta$ the local point of $Q$ on $k_*\hat H$ determined by $f\,,$ by
$\varphi$ the point of~$k_*\hat C_{H}^{G}(Q,f)$ such that $(\delta,\varphi)$ is the
local point of $Q = P$ on $k_*\hat Gb$ determined by $e\,,$ and by  $\hat G^\varphi$ the
converse image in~$\hat G$ of the stabilizer $C_{(\delta,\varphi)}$ of~$(\delta,\varphi)$
in~$C\,.$ Then $b$ is a block of~$\hat G^\varphi$ and, if~{\rm (Q)} holds for $(b,\hat
G^\varphi)\,,$ it holds for $(b,\hat G)\,.$\/} 

\medskip
\noindent
{\bf Proof:} According to £3.9, the pair $(\delta,\varphi)$ has been identified indeed with a local point of $Q = P$ 
on $k_*\hat Gb\,,$ so that $e$ determines $\varphi\,;$ moreover, it follows from equality~£3.2.2 that $C$ acts on 
the set of points of $k_*\hat C_{H}^{G}(Q,f)$ and therefore it makes sense to consider the converse image 
$\hat G^\varphi$ of $C_{(\delta,\varphi)}$ in $\hat G\,.$

\smallskip
Since $C_{\hat G}(Q,f)\i \hat K\i \hat G^\varphi\,,$ $b$ is also a block of 
$\hat G^\varphi$ [10,~Proposition~15.10] and we have $\hat C_{H}^{G^\varphi}(Q,f)= \hat C_{H}^{G}(Q,f)\,,$
so that $\varphi$ is also a point of $k_*\hat C_{H}^{G^\varphi}(Q,f)$ and we have $\hat E_{G^\varphi}(Q_{(\delta,\varphi)}) \cong
\hat E_{G}(Q_{(\delta,\varphi)})\,;$ moreover, by the {\it Frattini
argument\/}, the stabilizer ${\rm Aut}(\hat G)_{(P,e)}$ of $(P,e)$ in ${\rm Aut}(\hat
G)_b$ {\it covers\/}  ${\rm Out}_{k^*}(\hat G)_b$ and therefore we get a canonical group
homomorphism
$${\rm Out}_{k^*}(\hat G)_b\too {\rm Out}_{k^*}({\hat G}^\varphi)_b
\eqno £4.8.1.$$
Now, setting $\F^\varphi = \F_{\!(b,\hat G^{^\varphi})}\,,$ it follows from Corollary~£3.12 that we have
$\O {\rm Out}(\hat G)_b\-$module isomorphisms
$$\G_k (\F,\widehat\aut_{\F^{^{\rm nc}}})\cong \G_k (\F^\varphi,
\widehat\aut_{(\F^\varphi)^{^{\rm nc}}})\qq \G_k (\hat G,b) \cong \G_k (\hat G^\varphi,b)
\eqno £4.8.2.$$
Thus, if there is an $\O {\rm Out}_{k^*}(\hat G^\varphi)_b\-$module isomorphism 
$$\G_k (\F^\varphi,\widehat\aut_{(\F^\varphi)^{^{\rm nc}}}) \cong \G_k (\hat G^\varphi,b)
\eqno £4.8.3\phantom{.}$$
then we get an $\O{\rm Out}_{k^*} (\hat G)_b\-$module isomorphism $\G_k
(\F,\widehat\aut_{\F^{^{\rm nc}}})
\cong  \G_k (\hat G,b)\,.$ We are done.

\bigskip
\noindent
{\bf Proposition~£4.9.} {\it With the notation above, assume that $\hat G/\hat H$ is a $p'\-$group, that $\hat G/\hat K$ 
is cyclic and that  $\hat G = \hat K\.N_{\hat G}(Q_{(\delta,\varphi)})\,.$ Let  $\hat x$ be an element of 
$N_{\hat G}(Q_{(\delta,\varphi)})$ such that the image of $\hat x$ in $\hat G/\hat K$ is a generator of this quotient, set 
$\hat G' = \hat H\.\langle \hat x\rangle$ and choose  a block $b'$ of~$\hat G'$~such that ${\rm Br}_Q(b')\not= 0\,.$
If {\rm (Q)} holds for $(b',\hat G')$ then it holds for $(b,\hat G)\,.$\/}
\medskip
\noindent
{\bf Proof:} Let  $\varphi'$ be a point of the   $k\-$algebra 
 $k_*\hat  C_{H}^{G'}(Q,f)$ such that  $Q_{(\delta,\varphi')}$ is a defect pointed group of~$b'\,;$ since the quotient $\hat G'/\hat H$ is cyclic, it is clear that 
$\hat x$ normalizes $Q_{(\delta,\varphi')}$ and therefore we have 
$$E_{G'} (Q_{(\delta,\varphi')}) =  E_{G} (Q_{(\delta,\varphi)})
\eqno £4.9.1;$$
 moreover, since the quotient $E_{G} (Q_{(\delta,\varphi)})\big/\big(E_H(Q_\delta)\cap E_{G} (Q_{(\delta,\varphi)})\big)$
 is cyclic and, according to Theorem~£3.10 above,  the converse images of the intersection $E_H (Q_\delta)\cap E_{G} 
 (Q_{(\delta,\varphi)})$ in $\hat E_{G'} (Q_{(\delta,\varphi')})$ and  $\hat E_{G} (Q_{(\delta,\varphi)})$ admit a $k^*\-$group isomorphism lifting the identity,  it follows from Lemma~£4.10 below applied to the $k^*\-$extensions that equality~£4.9.1 
 also can be lifted to a $k^*\-$group isomorphism $\hat E_{G'} (Q_{(\delta,\varphi')})\cong \hat E_{G} (Q_{(\delta,\varphi)})\,.$ Consequently, it follows from Corollary~£3.12 that we have 
 canonical $\O\-$module isomorphisms
 $$\G_k (\hat G,b) \cong \G_k (\hat G',b')\qq \G_k (\F,\widehat\aut_{\F^{^{\rm nc}}})\cong 
 \G_k (\F',\widehat\aut_{\F'^{^{\rm nc}}})
\eqno £4.9.2.$$

\smallskip
Then, if {\rm (Q)} holds for $(b',\hat G')\,,$ it is clear that, for any block $(b'',\hat G'')$ isomorphic to $(b',\hat G')\,,$
we can choose an $\O\-$module isomorphism
$$\gamma_{(b'',\hat G'')} : \G_k (\hat G'',b'')\cong \G_k (\F'',\widehat\aut_{\F''^{^{\rm nc}}})
\eqno £4.9.3,$$
where $\F'' = \F_{\!(b'',\hat G'')}\,,$ in such a way that these isomorphisms are compatible with the isomorphisms between
these blocks. At this point, it is easily checked that the obvious composition
$$\G_k (\hat G,b)\cong \G_k (\hat G',b')\buildrel \gamma_{(b',\hat G')}\over\cong \G_k (\F',\widehat\aut_{\F'^{^{\rm nc}}})\cong \G_k (\F,\widehat\aut_{\F^{^{\rm nc}}})
\eqno £4.9.4\phantom{.}$$
is an $\O {\rm Out}_{k^*}(\hat G)_b\-$module isomorphism. We are done.

\bigskip
\noindent
{\bf Lemma~£4.10.} {\it  Let $K$ be a finite group, $H$ a normal subgroup of $K$ such that the quotient $K/H$
 is cyclic, $A$ a divisible Abelian group and $\hat H$ a central $A\-$extension of $H\,.$  Assume that the action of $K$ on $H$ 
 can be lifted to an action of~$K$ on $\hat H$ such that we have $\,\Bbb H^2(K/H,A) = \{0\}\,.$ Then, there exists an essentially unique $A\-$extension $\hat K$
of~$K$ containing $\hat H$ and lifting the inclusion map $H\to K\,.$ In particular, any
automorphism $\tau$ of $K$ stabilizing $H$ which can be lifted to an automorphism $\hat \sigma$ of $\hat H\,,$ can be lifted to an automorphism of $\hat K$ extending $\hat\sigma\,.$\/}

\medskip
\noindent
{\bf Proof:} Choose a cyclic subgroup $C$ of $K$ such that $K = H\.C$ and set 
$D = C\cap H\,;$ since the converse image $\hat D$ of $D$ in $\hat H$ is split, we can choose
 a splitting $\theta\,\colon D\to \hat D\i\hat H$ and, since $C\i K$ acts on $\hat H\,,$ we can consider
  the semidirect product $\hat H\!\rtimes C\,;$ inside, we define the ``inverse diagonal'' 
 $$\Delta^*(D) = \{(\theta (y),y^{-1})\}_{y\in D}
 \eqno £4.10.1\phantom{.}$$
and it is easily checked that $\Delta^*(D) $ is a subgroup  contained in the center of~$\hat H\!\rtimes C\,;$ 
then, it suffices to set 
$$\hat K = (\hat H\!\rtimes C)/\Delta^*(D)
\eqno £4.10.2;$$
indeed, the structural homomorphism $\hat H\to \hat H\!\rtimes C$ determines an injection $\hat H\to \hat K$
lifting the inclusion $H\i K\,.$

\medskip
Moreover, if $\hat{\phantom{K}}\!\! \!K$ is an $A\-$extension of~$K$ containing 
$\hat H$ and  lifting the inclusion map $H\to K\,,$ then $\hat{\phantom{K}}\!\! \!K * 
\hat K^\circ$ contains $\hat H *\hat H^\circ$ which is {\it canonically\/} isomorphic to
$A\times H$ and therefore, up to suitable identifications, the quotient
$(\hat{\phantom{K}}\!\! \!K * \hat K^\circ)/H$ is an $A\-$extension of the cyclic
group $K/H$ {\it via\/} the action which is induced by the action of $K$ on $\hat H\,;$ but,
we assume that $\,\Bbb H^2(K/H,A) = \{0\}\,;$ hence, this extension is split and therefore
 $\hat{\phantom{K}}\!\! \!K * \hat K^\circ$ is also split
or, more precisely, there is an isomorphism $\hat{\phantom{K}}\!\! \!K\cong \hat K$ inducing
the identity on $\hat H\,.$
\eject

\smallskip
In particular, if $\tau$ is  an automorphism  of $K$ stabilizing $H$ which can be lifted to an automorphism $\hat \sigma$ of $\hat H\,,$ it induces a group isomorphism 
$$\hat H\!\rtimes C\cong \hat H\!\rtimes \tau (C)
\eqno £4.10.3\phantom{.}$$
 mapping
$\Delta^*(D)$ onto $\Delta^*\big(\tau (D)\big)$ and therefore it determines an
isomorphism
$$\hat K\cong \big(\hat H\!\rtimes \tau (C)\big)\big/\Delta^*\big(\tau (D)\big)
\eqno £4.10.4;$$
but, the right member of this isomorphism is also an $A\-$extension 
of~$\hat H$ lifting the inclusion map $H\to K$ and therefore it admits an isomorphism to
$\hat K$ inducing the identity on $\hat H\,.$ We are done.

\medskip
£4.11. Thus,  we may assume that $\hat G/\hat H$ is a cyclic $p'\-$group. But, we have $\hat H\cong \hat{\prod}_{i\in I} \hat H_i$
and we want to reduce our situation to the case where $I$ has a unique element. In order to do this reduction, we will apply
[10,~Corollary~15.47] which forces us to move to a ``bigger'' situation; namely, for any $i\in I\,,$ let us denote by $K_i$ the image of
$\hat K$ in ${\rm Aut}(\hat H_i)\,;$ note that, by the very definition of $\hat K$ (cf.~£4.6), we have
$$K_i = H_i\.C_{K_i}(Q_i,f_i)
\eqno £4.11.1.$$
Since $K_i/H_i$ is cyclic, it follows from Lemma~£4.10 that there exists an essentially unique $k^*\-$group $\hat K_i$ containing $\hat H_i\,;$ set $\hat K^* = \hat{\prod}_{i\in I} \hat K_i\,.$ Then, since $K/H$ is cyclic, identifying $K$ to
its canonical image in $K^* = \prod_{i\in I}K_i\,,$  it follows again from Lemma~£4.10 that
we can identify $\hat K$ with the converse image of~$K$ in $\hat K^*\,.$

\medskip
£4.12.  Similarly, we can identify $G$ and $K^*$ with their 
image in ${\rm Aut}(\hat H)$ and, in this group, we set $G^* = K^*\.G\,;$ once again, since $G^*/K^*$ is cyclic, there exists an essentially unique $k^*\-$group $\hat G^*$ containing $\hat K^*$ and, since $\hat G/\hat H$ is cyclic, we can identify $\hat G$ with the converse image of~$G$ in $\hat G^*\,;$ then, it is clear that
$$\hat K^*\cap \hat G = \hat K\qq \hat K^* = \hat H\.C_{\hat G^*}(Q,f)
\eqno £4.12.1.$$
Moreover, for any $i\in I\,,$ since the quotient 
$$C_{H_i}^{K_i}(Q_i,f_i) = C_{K_i}(Q_i,f_i)/C_{H_i}(Q_i,f_i)
\eqno £4.12.2$$
is cyclic, the $k^*\-$group $\hat C_{H_i}^{K_i}(Q_i,f_i)$ is split and it is quite clear that we can choose a $N_G(Q,f)\-$stable family
of  $k^*\-$group homomorphisms
$$\hat\varphi_i : \hat C_{H_i}^{K_i}(Q_i,f_i)\too k^*
\eqno £4.12.3;$$
now, since $\hat C_H^{K^*} (Q,f)= \hat C_H^{G^*} (Q,f)\,,$ this family determines a $N_G(Q,f)\-$stable point $\varphi^*$
 of $k_*\hat C_H^{G^*} (Q,f)$ and then the pair $(\delta,\varphi^*)$ 
determines a local point\break
\eject
\noindent
 of~$Q$ on~$k_*\hat G^*$ (cf.~£3.9);  it is quite clear that  $Q_{(\delta,\varphi^*)}$ is a defect
pointed group of a block $b^*$ of $\hat G^*$ (cf.~£2.9) and we set $\F^* = \F_{\!(b^*,G^*)}\,.$
Now, we replace Proposition~16.25 in~[10] by the following result.

\bigskip
\noindent
{\bf Proposition~£4.13.} {\it With the notation above, assume that $\hat G/\hat H$ is a cyclic $p'\-$group and that  $\hat G = \hat K\.N_{\hat G}(Q_{(\delta,\varphi)})\,.$  If {\rm (Q)} holds for $(b^\star,\hat G^\star)$
then it holds for $(b,\hat G)\,.$\/}

\medskip
\noindent
{\bf Proof:} Since $N_G(Q,f)$ normalizes $Q_{(\delta,\varphi^*)}\,,$ we clearly have
$$E_{G^*} (Q_{(\delta,\varphi^*)}) =  E_{G} (Q_{(\delta,\varphi)})
\eqno £4.13.1;$$
 moreover, since the quotient $E_{G} (Q_{(\delta,\varphi)})\big/\big(E_H(Q_\delta)\cap E_{G} (Q_{(\delta,\varphi)})\big)$
 is cyclic and, according to Theorem~£3.10 above,  the converse images of the intersection $E_H (Q_\delta)\cap E_{G} 
 (Q_{(\delta,\varphi)})$ in $\hat E_{G^*} (Q_{(\delta,\varphi^*)})$ and  $\hat E_{G} (Q_{(\delta,\varphi)})$ admit a $k^*\-$group isomorphism lifting the identity,  it follows from Lemma~£4.10 that equality~£4.13.1 also can be lifted to 
a $k^*\-$group isomorphism $\hat E_{G^*} (Q_{(\delta,\varphi^*)})\cong \hat E_{G} (Q_{(\delta,\varphi)})\,.$ Consequently, it follows from Corollary~£3.12 that we have 
 canonical $\O\-$module isomorphisms
 $$\G_k (\hat G,b) \cong \G_k (\hat G^*,b^*)\qq \G_k (\F,\widehat\aut_{\F^{^{\rm nc}}})\cong 
 \G_k (\F^*,\widehat\aut_{(\F^*)^{^{\rm nc}}})
\eqno £4.13.2.$$

\smallskip
Then, if {\rm (Q)} holds for $(b^*,\hat G^*)\,,$ it is clear that, for any block $(\bar b^*,\hat G^*)$ isomorphic to $(b^*,\hat G^*)\,,$
we can choose an $\O\-$module isomorphism
$$\gamma_{(\bar b^*,\hat G^*)} : \G_k (\hat G^*,\bar b^*)\cong \G_k (\bar\F^*,\widehat\aut_{(\bar\F^*)^{^{\rm nc}}})
\eqno £4.13.3,$$
where $\bar\F^* = \F_{\!(\bar b^*,\hat G^*)}\,,$ in such a way that these isomorphisms are compatible with the isomorphisms between
these blocks. At this point, it is easily checked that the obvious composition
$$\G_k (\hat G,b)\cong \G_k (\hat G^*,b^*)\buildrel \gamma_{(b^*,\hat G^*)}\over
\cong \G_k (\F^*,\widehat\aut_{(\F^*)^{^{\rm nc}}})\cong \G_k (\F,\widehat\aut_{\F^{^{\rm nc}}})
\eqno £4.13.4\phantom{.}$$
is an $\O {\rm Out}_{k^*}(\hat G)_b\-$module isomorphism. We are done.

\medskip
£4.14. Consequently, from now on we assume that $K = \prod_{i\in I} K_i$ and that the quotients $C = G/K$ and $K_i/H_i$ 
for any $i\in I$ are cyclic $p'\-$groups; then, we have 
$$\hat K = \hat{\prod}_{i\in I}\hat K_i
\eqno £4.14.1\phantom{.}$$
and, since $b$ is also a block of $\hat K$ [10,~Proposition~15.10], we have $b = \otimes_{i\in I} b_i$ for a suitable
block $b_i$ of $\hat K_i$ for any $i\in I\,;$ moreover, we set $\K = \F_{\!(b,\hat K)}$ and $\K^i = \F_{\!(b_i,\hat K_i)}$
 for any $i\in I\,.$ Since $\hat G = \hat K\.N_{\hat G}(Q,f)\,,$ it is clear that
 $$C\cong \hat G/\hat K\cong N_{\hat G}(Q,f)/N_{\hat K}(Q,f)\cong \F (Q)/\K (Q)
 \eqno £4.14.2;$$
 then, for any subgroup $D$ of $C\,,$ we denote by ${}^D\!\hat K$ the converse image of $D$ in $\hat G$ and set
 ${}^D\K = \F_{\!(b,{}^D\!\hat K)}\,.$ Recall that we respectively denote by  $\R_{_{\!\!\hat K}}\!\G_k (\hat G,b)$
 and $\R_{_{\K}}\!\G_k(\F,\widehat\aut_{\F^{^{\rm nc}}})$ the intersection of
the kernels of all the respective $\O\-$module homomorphisms determined by the restriction
$$\G_k (\hat G,b)\to \G_k ({}^D\!\hat K,b)\qq \G_k(\F,\widehat\aut_{\F^{^{\rm nc}}})\to 
\G_k({}^D\K,\widehat\aut_{({}^D\K)^{^{\rm nc}}})
\eqno £4.14.3\phantom{.}$$
where $D$ runs over the set of proper subgroups of $C\,.$

 \medskip
 £4.15. It is clear that the quotient $C = G/K$ acts on $I\,;$ if $I$ decomposes on a disjoint union of
two nonempty $C\-$stable subsets $I'$ and $I''$ then, setting
$$\hat K' = \hat{\prod_{i'\in I'}} \hat K_{i'}\qq  \hat K'' = 
\hat{\prod_{i''\in I''}} \hat K_{i''}
\eqno £4.15.1,$$
  it follows again from Lemma~£4.10 that there exist essentially unique $k^*\-$groups 
$\hat G'$ and $\hat G''\,,$ respectively containing and normalizing $\hat K'$
and $\hat K''\,,$ such that
$$\hat G'/\hat K'\cong C\cong \hat G''/\hat K''
\qq \hat G'\,\hat\times_{C}\, \hat G''\cong \hat G
\eqno £4.15.2.$$
Moreover, setting $b' = \otimes_{i'\in I'} b_{i'}$ and $b'' = \otimes_{i''\in I''}
b_{i''}\,,$ it follows from [10~Pro-position~15.10] that $b'$ and $b''$ are respective
blocks of $\hat G'$ and $\hat G''\,;$ we set
$$\eqalign{\F' = \F_{\!(b',\hat G')} &\qq \F'' = \F_{\!(b'',\hat G'')}\cr
\K' = \F_{\!(b',\hat K')} &\qq \K'' = \F_{\!(b'',\hat K'')}\cr}
\eqno £4.15.3.$$
Note that, for any subgroup $D$ of $C\,,$ we have an analogous situation with
respect to the converse images ${}^D\!\hat K\,,$ ${}^{D}\!\hat K'$ and 
${}^{D}\!\hat K''$ of $D$ in $\hat G\,,$ $\hat G'$ and~$\hat G''\,.$

\bigskip
\noindent
{\bf Proposition~£4.16}\phantom{.} {\it  With the notation above, assume that  ${\rm Aut}_{k^*}(\hat G)_b$ stabilizes $I'$
and $I''\,.$ If {\rm (Q)} holds for $(b',{}^D\!\hat K')$ and $(b'',{}^D\!\hat K'')$
 for any subgroup $D$ of~$C$, then it holds for $(b,\hat G)\,.$\/}
\medskip
\noindent
{\bf Proof:} According to our hypothesis, we have canonical group homomorphisms
$${\rm Out}_{k^*}(\hat G')_{b'}\longleftarrow {\rm Out}_{k^*}(\hat G)_b\too
{\rm Out}_{k^*}(\hat G'')_{b''}
\eqno £4.16.1;$$
then, since any homomorphism from $C$ to $k^*$ induces $k^*\-$group
automorphisms of $\hat G\,,$ $\hat G'$ and $\hat G''$ which are contained in the
centers of
${\rm Aut}_{k^*}(\hat G)\,,$ ${\rm Aut}_{k^*}(\hat G')$ and~${\rm Aut}_{k^*}(\hat G'')$
respectively, it follows from [10,~Corollary~£15.47] that, setting $\hat\R = \R\G_k (C)\,,$ we have $\O {\rm Out}_{k^*}(\hat
G)_b\-$module isomorphisms
$$\eqalign{\R_{_{\!\hat K'}}\! \G_k (\hat G',b') \otimes_{\hat\R} 
\R_{_{\!\hat K''}}\!\G_k (\hat G''&,b'') \cong \R_{_{\!\hat K}}\!\G_k (\hat G,b)\cr 
\R_{_{\K'}}\!\G_k (\F',\widehat\aut_{\F'^{^{\rm nc}}})\otimes_{\hat\R} 
\R_{_{\K''}}\!\G_k (\F''&,\widehat\aut_{\F''^{^{\rm nc}}}) \cong \R_{_{\K}}\!\G_k (\F,\widehat\aut_{\F^{^{\rm nc}}})\cr}
\eqno £4.16.2.$$

\smallskip
Moreover, assume that we have $\O {\rm Out}_{k^*}(\hat G')_{b'}\-$ and $\O {\rm Out}_{k^*}(\hat G'')_{b''}\-$mo-dule isomorphisms
$$\G_k (\hat G',b')\cong \G_k (\F',\widehat\aut_{\F'^{^{\rm nc}}})\qq
\G_k (\hat G'',b'')
\cong \G_k (\F'',\widehat\aut_{\F''^{^{\rm nc}}}) 
\eqno £4.16.3;$$
since the restriction induces {\it compatible\/} $\G_k (C)\-$module structures on all the
members of these isomorphisms [10,~15.21 and~15.33], it follows from
[10, 15.23.2 and~£15.37.1] that we still have $\hat\R {\rm
Out}_{k^*}(\hat G)_b\-$module isomorphisms 
$$\eqalign{\R_{_{\!\hat K'}}\!\G_k (\hat G',b')&\cong  \R_{_{\!\K'}}\!\G_k (\F',\widehat\aut_{\F'^{^{\rm nc}}})\cr
\R_{_{\!\hat K''}}\!\G_k (\hat G'',b'')&\cong \R_{_{\!\K''}}\!\G_k (\F'',\widehat\aut_{\F''^{^{\rm nc}}})\cr}
\eqno £4.16.4.$$
Then, from isomorphisms~£4.16.2 we get an $\O {\rm Out}_{k^*}(\hat G)_b\-$module
isomorphism
$$\R_{_{\!\hat K}}\!\G_k (\hat G,b)\cong \R_{_{\K}}\!\G_k (\F,\widehat\aut_{\F^{^{\rm nc}}})
\eqno £4.16.5.$$

\smallskip
Consequently, according to our hypothesis, for any subgroup $D$ of $C$ we have
an  $\O {\rm Out}_{k^*}({}^{D}\!\hat K)_b\-$module isomorphism
$$\R_{_{\!\hat K}}\!\G_k ({}^{D}\!\hat K,b)\cong \R_{_{\K}}\!\G_k ({}^D \K,\widehat\aut_{({}^D \K)^{^{\rm nc}}}) 
\eqno £4.16.6;$$
but, since ${\rm Aut}_{k^*}(\hat G)_b$ stabilizes $\hat K\,,$ we have evident group
homomorphisms
$$C\too {\rm Out}_{k^*}({}^{\bar D}\!\hat K)_b\longleftarrow 
{\rm Aut}_{k^*}(\hat G)_b
\eqno £4.16.7\phantom{.}$$
and it is clear that the image of ${\rm Aut}_{k^*}(\hat G)_b$ contains and normalizes the
image of $C\,;$ hence, we still have an  $\O {\rm Out}_{k^*}(\hat G)_b\-$module
isomorphism
$$\R_{_{\!\hat K}}\!\G_k ({}^{D}\!\hat K,b)^{C}\cong\R_{_{\K}}\!\G_k ({}^D \K,\widehat\aut_{({}^D \K)^{^{\rm nc}}})^{C} 
\eqno £4.16.8.$$
Then, it follows from [10,~15.23.4 and~15.38.1] that the direct sum of
isomorphisms~£4.16.8 when $D$ runs over the set of subgroups of $C$ supplies an
 $\O {\rm Out}_{k^*}(\hat G)_b\-$module isomorphism 
$\G_k (\hat G,b)\cong \G_k (\F,\widehat\aut_{\F^{^{\rm nc}}}) \,.$ We are done.

\medskip
£4.17\phantom{.} From now on, we assume that the group ${\rm Aut}_{k^*}(\hat G)_b$ acts transitively on $I\,;$ in
particular, it acts transitively on the set of $C\-$orbits of $I$ and, for any $C\-$orbit~$O$ we consider the 
$k^*\-$group and the block
$$ \hat K^O = \hat{\prod}_{i\in O} \hat K_i\qq b^O =
\bigotimes_{i\in O} b_i
\eqno £4.17.1;$$
once again, it follows from Lemma~£4.10 that there exists an essentially unique $k^*\-$group
$\hat G^O$ containing $\hat K^O$ and fulfilling $C\cong \hat G^O/\hat K^O\,;$ then,
it follows from [10,~Proposition~15.10] that $b^O$ is also a block of $\hat G^O$ and we set
$$\F^O = \F_{\!(b^O,\hat G^O)}\qq \K^O = \F_{\!(b^O,\hat K^O)}
\eqno £4.17.2.$$ 
Note that $\hat G$ is isomorphic to the {\it direct sum over $C$\/} of the family of
$k^*\-$groups $\hat G^O$ when~$O$ runs over the set of $C\-$orbits of $I\,.$

\bigskip
\noindent
{\bf Proposition~£4.18}\phantom{.} {\it  With the notation above, assume that ${\rm Aut}_{k^*}(\hat G)_b$ acts transitively on
$I$ and let $O$ be a $C\-$orbit of~$I\,.$ If {\rm (Q)} holds for $\big(b^O,{}^{D}\!(\hat K^O)\big)$ for any subgroup $D$ of $C$,
 then it holds  for~$(b,\hat G)\,.$\/}
\medskip
\noindent
{\bf Proof:} It is clear that the action of ${\rm Aut}_{k^*}(\hat G)_b$ on $I$ induces an
action of ${\rm Out}_{k^*}(\hat G)_b$ on the set $\bar I$ of $C\-$orbits of $I\,;$
moreover, denoting by  ${\rm Out}_{k^*}(\hat G)_{b,O}$ the stabilizer of~$O$ in 
${\rm Out}_{k^*}(\hat G)_b\,,$ it is quite clear that the restriction induces a 
group homomorphism
$${\rm Out}_{k^*}(\hat G)_{b,O}\too {\rm Out}_{k^*}(\hat G^O)_{b^O}
\eqno £4.18.1\,.$$
On the other hand, it is quite clear that ${\rm Out}_{k^*}  (\hat G)_{b}$ acts transitively
on the two families of $\O\-$modules
$$\{\G_k  (\hat G^{O'},b^{O'})\}_{O'\in \bar I}\qq
\{\R_{_{\K^{O'}}}\!\G_k (\F^{O'},\widehat\aut_{(\F^{O'})^{^{\rm nc}}})\}_{O'\in \bar I}
\eqno £4.18.2\phantom{.}$$
and then, iterating the canonical isomorphisms in [10,~Corollary~15.47] and setting $\hat\R = \R\G_k(C)\,,$ it is not difficult to check that we have  $\hat\R {\rm Out}_{k^*}(\hat G)_b\-$mo-dule isomorphisms
$$\eqalign{{}^{\hat\R}{\bf Ten}\,_{{\rm Out}_{k^*} (\hat G)_{b,O}}^{{\rm Out}_{k^*} (\hat G)_b} \big(\R_{_{\!\hat K^O}}\! \G_k  (\hat G^O,b^O)\big) &\cong  \R_{_{\!\hat K}}\!\G_k (\hat G,b)\cr 
{}^{\hat\R}{\bf Ten}\,_{{\rm Out}_{k^*}(\hat G)_{b,O}}^{{\rm Out}_{k^*}
(\hat G)_b} \big(\R_{_{\K^O}}\!\G_k (\F^O,\widehat\aut_{(\F^O)^{^{\rm nc}}})\big)
&\cong \R_{_{\!\K}}\!\G_k (\F,\widehat\aut_{\F^{^{\rm nc}}})\cr}
\eqno £4.18.3$$
where ${}^{\hat\R}{\bf Ten}\,$ denotes the usual {\it tensor induction of $\hat\R\-$modules\/}.

\smallskip
Moreover, assume that we have an $\O {\rm Out}_{k^*}(\hat G^O)_{b^O}\-$module
isomorphism
$$\G_k (\hat G^O,b^O)\cong \G_k (\F^O,\widehat\aut_{(\F^O)^{^{\rm nc}}}) 
\eqno £4.18.4;$$
then, since the restriction induces compatible $\G_k (C)\-$module structures on both
members of this isomorphism [10,~15.21 and~15.33], it follows from
[10,~15.23.2 and~15.37.1] that we still have an $\O {\rm Out}_{k^*}(\hat G)_{b,O}\-$module isomorphism 
$$\R_{_{\!\hat K^O}}\!\G_k (\hat G^O,b^O)\cong \R_{_{\K^O}}\!\G_k (\F^O,\widehat\aut_{(\F^O)^{^{\rm nc}}}) 
\eqno £4.18.5.$$
Thus, from isomorphisms~£4.18.3 above, we get an $\O {\rm Out}_{k^*}(\hat
G)_b\-$module isomorphism
$$\R_{_{\!\hat K}}\!\G_k (\hat G,b)\cong \R_{_{\!\K}}\!\G_k (\F,\widehat\aut_{\F^{^{\rm nc}}})
\eqno £4.18.6.$$

\smallskip
Consequently, according to our hypothesis and possibly applying Proposition~£4.16 above and [10,~15.23.2 and~15.37.1], 
 for any subgroup $D$ of $C$ we have an  $\O {\rm Out}_{k^*} ({}^{D}\hat K)_b\-$module isomorphism
$$\R_{_{\!\hat K}}\!\G_k \big({}^{D}\!\hat K,b\big)\cong \R_{_{\!\K}}\!\G_k  ({}^D \K,\widehat\aut_{({}^D \K)^{^{\rm nc}}}\big) 
\eqno £4.18.7;$$
but, since ${\rm Aut}_{k^*}(\hat G)_b$ stabilizes $\hat K\,,$ we have evident group
homomorphisms
$$\bar C\too {\rm Out}_{k^*}({}^{D}\!\hat K)_b\longleftarrow 
{\rm Aut}_{k^*}(\hat G)_b
\eqno £4.18.8\phantom{.}$$
and it is clear that the image of ${\rm Aut}_{k^*}(\hat G)_b$ contains and normalizes the
image of $C\,;$ hence, we still have an  $\O {\rm Out}_{k^*}(\hat G)_b\-$module 
isomorphism
$$\R_{_{\!\hat K}}\!\G_k \big({}^{D}\!\hat K,b\big)^C \cong 
\R_{_{\!\K}}\!\G_k  ({}^D \K,\widehat\aut_{({}^D \K)^{^{\rm nc}}}\big)^C
\eqno £4.18.9.$$
Then, it follows from [10,~15.23.4 and~15.38.1] that the direct sum of
isomorphisms~£4.18.9 when $D$ runs over the set of subgroups of $C$ supplies 
an $\O {\rm Out}_{k^*}(\hat G)_b\-$module isomorphism 
$\G_k (\hat G,b)\cong \G_k (\F,\widehat\aut_{\F^{^{\rm nc}}})\,.$ We are done.

\medskip
£4.19.  In the last step of our reduction, we  assume that $C = G/K$ acts transitively on~$I\,.$ 
 In this situation, we have to consider the {\it direct product\/} of groups
$$\hat{\hat K} = \prod_{i\in I} \hat K_i
\eqno £4.19.1;$$
 since $C$ is cyclic and it acts on $(k^*)^I$ permuting the factors, it follows from Lemma~£4.10 
 that there exists an essentially unique  $(k^*)^I\-$extension $\skew3\hat{\hat G}$ of $G$ containing 
 $\hat{\hat K}\,.$
Moreover, denoting by  $\nabla_{\!k^*}\,\colon (k^*)^I\to k^*$  the group homomorphism
induced by the product in $k$ and   considering the {\it group algebras\/} of the {\it groups\/} $(k^*)^I$ and
$\skew3\hat{\hat G}$ over $k$ and the
$k\-$algebra homomorphism $k(k^*)^I\to k$ determined by~$\nabla_{\!k^*}\,,$ it is quite
clear that we have a $k^*\-$group and a $k\-$algebra isomorphisms
$$\skew3\hat{\hat G}/{\rm Ker}(\nabla_{\!k^*})\cong \hat G\qq k\otimes_{k(k^*)^I}
k\skew3\hat{\hat G}\cong k_*\hat G
\eqno £4.19.2.$$

\par
£4.20. In particular, since $\hat G$ acts transitively on the
family  $\{\hat K_i\}_{i\in I}$ and, for any $i\in I\,,$ $b_i$ is a block of $\hat K_i$
(cf.~£4.14), by the {\it Frattini argument\/} we get canonical group homomorphisms
$${\rm Out}_{k^*}(\hat G)_b\too {\rm Out}_{k^*}(\hat K_i)_{b_i}
\eqno £4.20.1.$$
Moreover, choose an element $i\in I$ and respectively denote by $C_i\,,$ $\hat G_i$ and $\skew3\hat{\hat G}_i$ 
the stabilizers of $i$ in $C\,,$ $\hat G$ and $\skew3\hat{\hat G}\,,$ which actually act trivially on $I\,;$ setting $I' = I -\{i\}\,,$ it is clear that  $\prod_{i'\in I'} \hat K_{i'}$ is a normal subgroup of $\skew3\hat{\hat G}_i$ and that the quotient
$$\hat G^i = \skew3\hat{\hat G}_i \big/\big(\prod_{i'\in I'} \hat K_{i'}\big)
\eqno £4.20.2\phantom{.}$$
is a $k^*\-$group which contains $\hat K_i$ as a normal $k^*\-$subgroup;  now, any
$k_*\hat G^i\-$mo-dule $M_i$ can be viewed as a  $k \skew3\hat{\hat G}_i \-$module, and the 
point is that {\it the tensor induction 
${\rm Ten}_{ \skew3\hat{\hat G}_i }^{\skew3\hat{\hat G}} (M_i)$ becomes a $k_*\hat G\-$module\/} Let us consider 
$\R\G_k (C)$ as an $\R\G_k (C_i)\-$algebra {\it via\/} the group homomorphism  mapping $c\in C$ on~$c^{\vert I\vert}\,.$

\bigskip
\noindent 
{\bf Proposition~£4.21.} {\it With the notation above, assume that $C$ acts transitively on $I$ and  choose an element $i$ 
of $I\,.$  For any   $k_*\hat G^i\-$module $M_i$  considered as a  $k\skew3\hat{\hat G}_i\-$module,    
${\rm Ten}_{\skew3\hat{\hat G}_i}^{\skew3\hat{\hat G}} (M_i)$  becomes a $k_*\hat G\-$module, and  
this correspondence induces an $\O{\rm Out}_{k^*}(\hat G)_b\-$module isomorphism
$$\R\G_k (C)\otimes_{\R\G_k (C_i)} \R_{_{\!\hat K^i\!\!}}\G_k  (\hat G^i,b_i)\cong \R_{_{\!\hat K\!}} \G_k  (\hat G,b)
\eqno £4.21.1.$$\/}

\par
\noindent
{\bf Proof:} Recall that the {\it tensor induction\/} of $M_i$ from~$\skew3\hat{\hat G}_i$
to~$\skew3\hat{\hat G}$ is the $k\skew3\hat{\hat G}\-$module
$${\rm Ten}_{\skew3\hat{\hat G}_i}^{\skew3\hat{\hat G}} (M_i) = \bigotimes_{X\in
\skew3\hat{\hat G}/\skew3\hat{\hat G}_i} (kX\otimes_{k\skew3\hat{\hat G}_i} M_i)
\eqno £4.21.2,$$
where $kX$ denotes the $k\-$vector space over the (right-hand) $\skew3\hat{\hat G}_i\-$class $X$ of~$\skew3\hat{\hat G}\,,$ endowed with the (right-hand) $k\skew3\hat{\hat G}_i\-$module structure
determined by the multiplication on the right  [10,~8.2]. It is clear~that in ${\rm Ten}_{\skew3\hat{\hat G}_i}^{\skew3\hat{\hat
G}} (M_i)$ the multiplication by $(\lambda_i)_{i\in I}\in (k^*)^I$ coincides with the
multiplication by
$\prod_{i\in I}\lambda_i\in k^*\,,$ so that, according to isomorphisms~£4.19.2, 
${\rm Ten}_{\skew3\hat{\hat G}_i}^{\skew3\hat{\hat G}} (M_i)$ becomes a $k_*\hat
G\-$module.

\smallskip
Moreover,  if we have $M_i = M'_i\oplus M''_i$ as $k\skew3\hat{\hat G}_i\-$modules then we clearly get
$$\eqalign{{\rm Ten}_{\skew3\hat{\hat G}_i}^{\skew3\hat{\hat G}} &(M_i)\cr
& = \bigoplus_{\X} \big(\bigotimes_{X\in \X}(kX\otimes_{k\skew3\hat{\hat G}_i} M'_i)\big)\otimes_k \big(\bigotimes_{Y\in \skew3\hat{\hat G}/\skew3\hat{\hat G}_i -\X}(kY\otimes_{k\skew3\hat{\hat G}_i} M''_i)\big)\cr}
\eqno~£4.21.3\phantom{.}$$
where $\X$ runs over the set of all the subsets of $\skew3\hat{\hat G}/\skew3\hat{\hat G}_i\,;$ in particular, 
since any $p'\-$element  $\skew1\hat{\hat x}\in  \skew3\hat{\hat G}$ such that $\skew3\hat{\hat G} = \skew3\hat{\hat G}_i\.\langle \skew1\hat{\hat x} \rangle$ stabilizes this direct sum but only fixes the terms labeled by $\emptyset$ and by
$\skew3\hat{\hat G}/\skew3\hat{\hat G}_i\,,$ denoting by $\chi_i$ the {\it modular character\/} of $M_i\,,$ we get
$$\big({\rm Ten}_{\skew3\hat{\hat G}_i}^{\skew3\hat{\hat G}} (\chi_i)\big)(\skew1\hat{\hat
x}) = \chi_i \big(\,\overline{\skew1\hat{\hat x}^{\vert I \vert}}^i\big)
\eqno £4.21.4\phantom{.}$$
\eject
\noindent
where $\,\overline{\skew1\hat{\hat x}^{\vert I \vert}}^i$ denotes the image
of  $\skew1\hat{\hat x}^{\vert I \vert}\in  \skew3\hat{\hat G}_i$ in $\hat G^i\,;$ in
particular, this equality shows that the tensor induction ${\rm Ten}_{\skew3\hat{\hat G}_i}^{\skew3\hat{\hat G}}$
 induces an $\O\-$module homomorphism from $\G_k (\hat G^i)$ to
the $\O\-$module formed by the restriction  of {\it modular characters\/} of $\hat G$ to the set of
$p'\-$elements $\hat x\in \hat G$ such that $\hat G = \hat G_i\.\langle\hat x\rangle\,.$

\smallskip
But, by the very definition of $\R_{_{\!\hat G_i\!}}\G_k  (\hat G)$ in [10,~15.22.4], $\R_{_{\!\hat G_i\!}}\G_k  (\hat G)$ is isomorphic to this $\O\-$module and this restriction is equivalent to the projection obtained from
[10,~£15.23.4]
$$\G_k  (\hat G)\too \R_{_{\!\hat G_i\!}}\G_k  (\hat G)
\eqno £4.21.5,$$
so that, we finally get an $\O\-$module homomorphism
$$\G_k  (\hat G^i)\too \R_{_{\!\hat G_i\!}}\G_k  (\hat G)
\eqno £4.21.6;$$
more precisely, it is quite clear that this homomorphism maps $\R_{_{\!\hat K^i\!}}\G_k  (\hat G^i)$
on~$\R_{_{\!\hat K\!}}\G_k  (\hat G)$ and, since $\R_{_{\!\hat K^i\!}}\G_k  (\hat G^i)$ and $\R_{_{\!\hat K\!}}\G_k  (\hat G)$ 
respectively have  $\R\G_k (C_i)\-$ and   $\R\G_k (C)\-$module structures, which are compatible with homomorphism £4.21.5
above, we still get an $\R\G_k (C)\-$module homomorphism
$$\R\G_k (C)\otimes_{\R\G_k (C_i)} \R_{_{\!\hat K^i\!}}\G_k  (\hat G^i)\too 
\R_{_{\!\hat K\!}}\G_k  (\hat G)
\eqno £4.21.7\phantom{.}$$
and we claim that it is bijective.

\smallskip
Indeed, if $M_i$ is a simple $k_*\hat G^i\-$module such that the restriction to $\hat K_i$ remains simple
 then, since $\hat K = \widehat{\prod}_{j\in I} \hat K_j$ and therefore $k_*\hat K\cong \otimes_{j\in I} k_*\hat K_j\,,$ 
 it is clear that the restriction of  ${\rm Ten}_{\skew3\hat{\hat G}_i}^{\skew3\hat{\hat G}} (M_i)$ to $k_*\hat K$ is
simple too. Conversely, if $M$ is a simple $k_*\hat G\-$module such that the restriction to $\hat K$ remains simple then we necessarily
have $M \cong \otimes_{j\in I} M_j$ or, more explicitly, 
$${\rm End}_k (M)\cong \bigotimes_{j\in I} S_j
\eqno £4.21.8\phantom{.}$$
where $S_j\cong {\rm End}_k (M_j)$ is generated by the image of $\hat K_j\i \hat K$ for any $j\in I\,;$ thus, 
$\hat G$ stabilizes the family of $k\-$subalgebras $\{S_j\}_{j\in I}$ of ${\rm End}_k (M)\,,$ and~$\hat G_j\,,$ which coincides with
$\hat G_i\,,$ stabilizes  $S_j$  for any $j\in I\,;$ in particular, the image of any $\hat x_i\in \hat G_i$ in ${\rm End}_k (M)$
has the form $\otimes_{j\in I} s_j$ for suitable $s_j\in S_j$   for any $j\in I\,;$ then, considering $\hat G_i$ as a quotient of 
$\skew3\hat{\hat G}_i\,,$ it is clear that the corresponding homomorphism $\skew3\hat{\hat G}_i\to {\rm End}_k (M)$ factorizes throughout a $(k^*)^I\-$extension homomorphism 
$$\skew3\hat{\hat G}_i\too \prod_{j\in I} S_j
\eqno £4.21.9\phantom{.}$$
and then that the corresponding homomorphism $\skew3\hat{\hat G}_i\to S_i$  factorizes throughout a $k^*\-$group
homomorphism $\hat G^i\to S_i\,,$ so that $M_i$ becomes a $k_*\hat G^i\-$module which remains simple restricted to $\hat K_i\,.$

\smallskip
Moreover,  the groups ${\rm Hom}(C_i,k^*)$ and ${\rm Hom}(C,k^*)$ respectively determine $k^*\-$group
automorphisms of $\hat G^i$ and   $\hat G\,,$ and it is elementary to check that the restriction of $M$ {\it via\/} the 
 $k^*\-$group automorphism of $\hat G$ determined by a suitable element of ${\rm Hom}(C,k^*)$ coincides with 
 ${\rm Ten}_{\skew3\hat{\hat G}_i}^{\skew3\hat{\hat G}} (M_i)\,.$ Now, it is not difficult to check that this
correspondence induces a bijection between the set of~${\rm Hom}(C_i,k^*)\-$orbits of isomorphism classes of simple
$k_*\hat G^i\-$modules which remain simple restricted to~$\hat K_i\,,$ and the set of ${\rm Hom}(C,k^*)\-$orbits of isomorphism
classes of the simple $k_*\hat G\-$modules which remain simple restricted to~$\hat K\,.$

\smallskip
But, it follows from [10,~15.23.2] that these sets of isomorphism classes respectively label $\R\G_k (C_i)\-$ and   
$\R\G_k (C)\-$bases of $\R_{_{\!\hat K^i\!}}\G_k  (\hat G^i)$ and $ \R_{_{\!\hat K\!}}\G_k  (\hat G)\,;$
 this implies the bijectivity of homomorphism~£4.21.6 above, proving our claim. Moreover, it is easily checked that $M_i$  is associated with the block   $b_i$ if and only if ${\rm Ten}_{\skew3\hat{\hat G}_i}^{\skew3\hat{\hat G}} (M_i)$ is associated
with the block~$b\,;$ hence, isomorphism~£4.21.6 induces the $\O\-$module
isomorphism~£4.21.1.

\smallskip
 On the other hand, since any $k^*\-$automorphism $\hat\sigma$ of $\hat G$ stabilizes the
family $\{\hat H_j\}_{j\in I}\,,$ if $\hat\sigma$ stabilizes $b$ then it also stabilizes
both families $\{\hat K_j\}_{j\in I}$ and $\{\hat G^j\}_{j\in I}\,,$ and therefore,
according to Lemma~£4.10, $\hat\sigma$ can be lifted to an automorphism $\hat{\hat\sigma}$
of $\skew3\hat{\hat G}\,;$ moreover, since $\hat G$ acts transitively on $I\,,$ up to a
modification of
$\hat\sigma$ by an inner automorphism of $\hat G\,,$ we may assume that $\hat\sigma$
fixes $i$ and then $\hat{\hat\sigma}$ determines a $k^*\-$automorphism $\hat\sigma_i$ of
$\hat G^i\,;$ in this case, assuming that $M_i$ is associated with the block $b_i\,,$ it
is quite clear that
$${\rm Ten}_{\skew3\hat{\hat G}_i}^{\skew3\hat{\hat G}}\big({\rm Res\,}_{\hat\sigma_i}(M_i)\big)
\cong {\rm Res}_{\,\hat{\hat\sigma}}\big({\rm Ten}_{\skew3\hat{\hat G}_i}^{\skew3\hat{\hat G}}
(M_i)\big)
\eqno £4.21.10\phantom{.}$$
and therefore, since homomorphism~£4.20.1 maps the class of $\hat\sigma_i$ on the class
of~$\hat\sigma\,,$ homomorphism~£4.21.6 is actually an $\O {\rm Out}_{k^*}(\hat
G)_b\-$module homomorphism. We are done.

\medskip
£4.22. Now, always choosing an element $i$ in $I$ and setting $\F^i = \F_{\!(b_i,\hat G^i)}\,,$ we have an analogous 
result on the relationship between $\R_{_{\K}}\! \G_k (\F,\widehat{\aut}_{\F^{^{\rm nc}}})$ and 
$ \R_{_{\K^i}}\!\G_k (\F^i,\widehat{\aut}_{(\F^i)^{^{\rm nc}}})\,;$  here we need the alternative definition~£2.12.3.

\bigskip
\noindent
{\bf Theorem~£4.23.} {\it  With the notation above, assume that $C$ acts transitively on $I$ and  choose an element $i$ 
of $I\,.$ Then,  there is an $\O {\rm Out}_{k^*}(\hat G)_b\-$module isomorphism
$$\R\G_k (C)\otimes_{\R\G_k (C_i)} \R_{_{\K^i}}\!\G_k (\F^i,\widehat{\aut}_{(\F^i)^{^{\rm nc}}})
\cong \R_{_{\K}}\! \G_k (\F,\widehat{\aut}_{\F^{^{\rm nc}}})
\eqno £4.23.1.$$\/}

\par
\noindent
{\bf Proof:} Let us recall our notation in [10,~15.33]; let $\frak r\,\colon \Delta_n\to \F^{^{\rm sc}}$  be a  
$\F^{^{\rm sc}}\-$chain  (cf.~£2.12); since $\frak r (n)$ is also $\K\-$selfcentralizing [10,~Lemma~15.16] and we
identify $\F (\frak r)$ with the stabilizer in $\F\big(\frak r (n)\big)$ of the images of $\frak r(\ell)$ in $\frak r (n)$ 
for any $\ell\in \Delta_n\,,$ it makes sense to consider $\K (\frak r)= \K\big(\frak r (n)\big)\cap \F (\frak r)$ which is a normal
subgroup of $\F (\frak r)\,;$ then, we denote by $\ch^*_C(\F^{^{\rm sc}}\big)$ the full subcategory of 
$\ch^* (\F^{^{\rm sc}})$ over the set of $\F^{^{\rm sc}}\-$chains $\frak r\,\colon \Delta_n\to \F^{^{\rm sc}}$ such that 
$$\F (\frak r)/\K (\frak r) \cong C
\eqno £4.23.2.$$ 
More explicitly, by the very definition of $\F = {}^C \K\,,$ we have 
$$\F (Q)/\K (Q)\cong C
\eqno £4.23.3;$$
choosing a lifting $\sigma \in \F(Q)$ of a generator of $C\,,$ we have a {\it Frobenius functor\/} 
$\frak f_\sigma\,\colon \F\to \F$ [10,~12.1] and therefore we get a new  $\F^{^{\rm sc}}\-$chain
 $\frak f_\sigma\circ \frak r\,;$ then, isomorphism~£4.23.2 is equivalent to the existence of a {\it natural\/} isomorphism
 $\nu\,\colon \frak r\cong \frak f_\sigma\circ \frak r$ formed by $\K\-$isomorphisms.

\smallskip
{\it Mutatis mutandis\/}, we also consider the corresponding full subcategory  $\ch^*_{C_i} \big((\F^i)^{^{\rm sc}}\big)$
of $\ch^* \big((\F^i)^{^{\rm sc}}\big)\,.$ Then, it follows from [10,~15.36] that we have   {\it contravariant\/} functors
$$\eqalign{\R_{_{\hat\K (\bullet)}}\! \G_k \big(\hat\F (\bullet)\big) &:
\ch^*_{C} (\F^{^{\rm sc}}\big)\too \O\-\mod\cr
\R_{_{\hat\K^i (\bullet)}}\! \G_k \big(\hat\F^i (\bullet)\big) &: 
\ch^*_{C_i}\big((\F^i)^{^{\rm sc}}\big)\too \O\-\mod\cr}
\eqno £4.23.4\phantom{.}$$
respectively mapping any $\ch^*_C (\F^{^{\rm sc}}\big)\-$object $(\frak r,\Delta_n)$ on 
$\R_{_{\hat\K (\frak r)}}\! \G_k \big(\hat\F (\frak r)\big)$ and any $\ch^*_{C_i} \big((\F^i)^{^{\rm sc}}
\big)\-$object  $(\frak r_i,\Delta_n)$ on  $\R_{_{\hat\K^i (\frak r_i)}}\! \G_k \big(\hat\F^i (\frak r_i)\big)\,,$
and from [10,~Proposi-tion~15.37] that we still have
$$\eqalign{\R_{_{\K}}\!\G_k(\F,\widehat \aut_{\F^{^{\rm nc}}})&\cong 
\lim_{\longleftarrow}\, \R_{_{\hat\K (\bullet)}}\! \G_k \big(\hat\F (\bullet)\big)\cr
\R_{_{\K^i}}\!\G_k(\F^i,\widehat \aut_{(\F^i)^{^{\rm nc}}}) &\cong 
\lim_{\longleftarrow}\, \R_{_{\hat\K^i (\bullet)}}\! \G_k \big(\hat\F^i (\bullet)\big)\cr}
\eqno £4.23.5.$$

\smallskip
Explicitly, we have
$$\R_{_{\K}}\!\G_k(\F,\widehat \aut_{\F^{^{\rm nc}}})\i \prod_{\frak r} \R_{_{\hat\K (\frak r)}}\! 
\G_k \big(\hat\F (\frak r)\big)
\eqno £4.23.6,$$
where $(\frak r,\Delta_m)$ runs over the set of $\ch^*_C(\F^{^{\rm sc}})\-$objects, and the left member coincides with
the set of $(X_\frak r)_\frak r\,,$ where $X_\frak r\in \R_{_{\hat\K (\frak r)}}\! \G_k \big(\hat\F (\frak r)\big)\,,$
which are ``stable'' by $\ch^*_C(\F^{^{\rm sc}}\big)\-$isomorphisms and, for any $\ch^*_C(\F^{^{\rm sc}}\big)\-$object
$(\frak q,\Delta_n)$ such that $\frak q = \frak r\circ \iota$ for some injective order-preserving map
$\iota\,\colon \Delta_n\to \Delta_m\,,$ the corresponding restriction map sends $X_\frak q$ to $X_\frak r\,;$ in particular,
 it is clear that we can~restrict ourselves to the $\ch^*_C(\F^{^{\rm sc}}\big)\-$objects such that we have 
 $\frak r (\ell -1)\i \frak r(\ell)$ and $\frak r(\ell\! -\!1\bullet\ell)$ is the inclusion map for any $1\le \ell \le m\,.$
 \eject

\smallskip
Moreover, for any $\ell\in \Delta_m\,,$ let us denote by $\frak  r_j(\ell)$ the image of $\frak r(\ell)$ in $Q_j$ and consider the
$\F^{^{\rm sc}}\-$chain $\frak r^*\,\colon \Delta_m\to  \F^{^{\rm sc}}$ mapping $\ell\in \Delta_m$ 
on $\prod_{j\in I} \frak r_j(\ell)$ and the $\Delta_m\-$morphisms on the corresponding inclusions; since $\K$ is {\it normal\/}
in~$\F$ [10,~12.6], it follows from [10,~Proposition~12.8] that any $\F\-$automorphism of $\frak r (m)$ induces
an $\F\-$automorphism of $\frak r^* (m)\,;$ in particular, we get a group homomorphism $\F (\frak r )\to \F (\frak r^*)$ 
and therefore $\frak r^*$ also fulfills condition~£4.23.2; furthermore, our functor $\widehat{\aut}_{\F^{^{\rm nc}}}$ lifts this homomorphism  to a $k^*\-$group homomorphism $\hat\F (\frak r )\to \hat\F (\frak r^*)\,.$

\smallskip
Consequently, considering the corresponding full subcategory, it follows from [10,~Proposition~A4.7] that, in the direct product
in~£4.23.6 above, we can restrict ourselves to the $\F^{^{\rm sc}}\-$chains $\frak r$ such that
we have $\frak r^* = \frak r\,;$ more explicitly, we may assume that, for any $\ell\in \Delta_m\,,$ we have
$$\frak r (\ell) = \prod_{j\in I} \frak r_j (\ell)
\eqno £4.23.7\phantom{.}$$
where, for any $j\in I\,,$ $\frak r_j\,\colon \Delta_m\to (\F^j)^{^{\rm sc}}$ is a 
$(\F^j)^{^{\rm sc}}\-$chain, fulfilling the corresponding condition~£4.23.2,  such that we have
 $\frak r_j (\ell -1)\i \frak r_j(\ell)$ and $\frak r_j(\ell\! -\!1\bullet\ell)$ is the inclusion map for any $1\le \ell \le m\,.$

 \smallskip
 In this case, it is quite clear that
 $$\hat \K (\frak r)\cong \hat{\prod}_{j\in I} \hat \K^j(\frak r_j)
 \eqno £4.23.8\phantom{.}$$
 and, arguing as in~£4.19 and~£4.20 above,   it follows from Proposition~£4.21 above that we have a {\it canonical\/} 
 $\R\G_k (C)\-$isomorphism
 $$\rho_\frak r : \R_{_{\hat\K (\frak r)}}\!  \G_k \big(\hat\F (\frak r)\big)\cong \R\G_k (C) \otimes_{\R\G_k (C_i)}
 \R_{_{\hat\K^i (\frak r_i)}}\!  \G_k \big(\hat\F^i (\frak r_i)\big)
\eqno £4.23.9\phantom{.}$$

\smallskip
{\it Mutatis mutandis\/}, we have
$$\R_{_{\K^i}}\!\G_k(\F^i,\widehat \aut_{(\F^i)^{^{\rm nc}}})\i \prod_{\frak r_i} \R_{_{\hat\K^i (\frak r_i)}}\! 
\G_k \big(\hat\F^i (\frak r_i)\big)
\eqno £4.23.10\phantom{.}$$
where $(\frak r_i,\Delta_m)$ runs over the set of $\ch^*_{C_i} \big((\F^i)^{^{\rm sc}}\big)\-$objects; once again,
 we can restrict ourselves to the $\ch^*_{C_i}\big((\F^i)^{^{\rm sc}}\big)\-$objects such that 
 $\frak r_i (\ell -1)\i \frak r_i(\ell)$ and $\frak r_i(\ell\! -\!1\bullet\ell)$ is the inclusion map for any $1\le \ell \le m\,.$ In particular,
 the extension 
 $$\R\G_k (C) \otimes_{\R\G_k (C_i)} \R_{_{\K^i}}\!\G_k(\F^i,\widehat \aut_{(\F^i)^{^{\rm nc}}})
 \eqno £4.23.11\phantom{.}$$
 coincides with the set of $(X_{\frak r_i})_{\frak r_i}\,,$ where $(\frak r_i,\Delta_m)$ runs over the set of such 
$\ch^*_{C_i}\big((\F^i)^{^{\rm sc}}\big)\-$objects and $X_{\frak r_i}$ belongs to the extension 
$$\R\G_k (C) \otimes_{\R\G_k (C_i)}\R_{_{\hat\K^i (\frak r_i)}}\! \G_k \big(\hat\F^i (\frak r_i)\big)
\eqno £4.23.12,$$
which are ``stable'' by $\ch^*_{C_i}\big((\F^i)^{^{\rm sc}}\big)\-$isomorphisms and, moreover, for such a
$\ch^*_{C_i}\big((\F^i)^{^{\rm sc}}\big)\-$object $(\frak q_i,\Delta_n)$ fulfilling $\frak q_i = \frak r_i\circ \iota$ for some injective order-preserving map $\iota\,\colon \Delta_n\to \Delta_m\,,$ the corresponding restriction map sends $X_{\frak q_i}$ to 
$X_{\frak r_i}\,.$

\smallskip
On the other hand, for any subgroup $R_i$ of $Q_i\i Q\,,$ we consider the subgroup of $Q$
$$R_i^* = \prod_{0\le t < \vert I\vert} \sigma^t (R_i)
\eqno £4.23.13\phantom{.}$$
and, more generally, for any $(\F^i)^{^{\rm sc}}\-$chain $\frak r_i$ defined by inclusion maps, we denote by $\frak r_i^*$ the corresponding $\F^{^{\rm sc}}\-$chain. Further, if $\frak r_i$ is a $\ch^*_{C_i}\big((\F^i)^{^{\rm sc}}\big)\-$ object
then $\frak r_i^*$ is a $\ch^*_{C} (\F^{^{\rm sc}})\-$object; indeed, the $\F^{^{\rm sc}}\-$chain
$\frak f_\sigma\circ \frak r_i^*$ maps $\ell\in\Delta_m$ on 
$$\sigma^{\vert I\vert}\big(\frak r_i(\ell)\big)\times  \prod_{1\le t < \vert I\vert} \sigma^t\big(\frak r_i(\ell)\big)
\eqno £4.23.14;$$
 but, we are assuming that there is a {\it natural\/} isomorphism
$\frak r_i\cong \frak f_{\sigma^{\vert I\vert}}\circ \frak r_i$ formed by $\K^i\-$isomorphisms; hence, we have
 a {\it natural\/} isomorphism $\frak r_i^*\cong \frak f_{\sigma}\circ \frak r_i^*$ formed by $\K\-$isomorphisms.

 \smallskip
 Moreover, it is quite clear that any $\ch^*\big((\F^i)^{^{\rm sc}}\big)\-$isomorphism $\frak r_i\cong\frak r'_i$ can be lifted to a 
 $\ch^*(\F^{^{\rm sc}})\-$isomorphism $\frak r_i^*\cong\frak r'^*_i\,,$ and that $\frak q_i = \frak r_i\circ \iota$
forces $\frak q^*_i = \frak r^*_i\circ \iota\,.$ Consequently, if $(X_\frak r)_\frak r$ is an element of 
$\R_{_{\K}}\!\G_k(\F,\widehat \aut_{\F^{^{\rm nc}}})$ then $\big(\rho_{\frak r_i^*}(X_{\frak r_i^*})\big)_{\frak r_i}$
is clearly an element of $\R_{_{\K^i}}\!\G_k(\F^i,\widehat \aut_{(\F^i)^{^{\rm nc}}})\,.$

\smallskip
Conversely, for any $\ch^*_{C} (\F^{^{\rm sc}})\-$object $(\frak r,\Delta_m)$ such that $\frak r$ is defined by inclusion maps and fulfills equality~£4.23.7, it is quite clear that the corresponding $(\F^i)^{^{\rm sc}}\-$chain $\frak r_i$ is also
defined by inclusion maps and that $(\frak r_i,\Delta_m)$ is a $\ch^*_{C_i}\big((\F^i)^{^{\rm sc}}\big)\-$object;
then, from isomorphism~£4.23.2 it is easy to check that $\frak r_i^*$ is {\it naturally\/} isomorphic to $\frak r$ {\it via\/} an isomorphism inducing the identity on $\frak r_i\,.$ Finally, for any element $(X_{\frak r_i})_{\frak r_i}$ of the extension~£4.23.12 
 and any $\ch^*_{C} (\F^{^{\rm sc}})\-$object $(\frak r,\Delta_m)$ as above, 
we can define $X_\frak r$ as the image of $(\rho_{\frak r_i^*})^{-1}(X_{\frak r_i})$ by the $\R\G_k (C)\-$module isomorphism
$$\R_{_{\hat\K (\frak r_i^*)}}\!  \G_k \big(\hat\F (\frak r_i^*)\big)\cong \R_{_{\hat\K (\frak r)}}\!  \G_k \big(\hat\F (\frak r)\big)
\eqno £4.23.15\phantom{.}$$ 
determined by a {\it natural\/} isomorphism $\frak r_i^*\cong \frak r$ inducing the identity on~$\frak r_i\,;$ it is easily checked that 
this definition does not depend on the choice of the {\it natural\/} isomorphism $\frak r_i^*\cong \frak r$ inducing the identity on 
$\frak r_i\,,$ and that the element $(X_\frak r)_\frak r$ of the corresponding direct product actually belongs to 
$\R_{_{\K}}\!\G_k(\F,\widehat \aut_{\F^{^{\rm nc}}})$ (cf.~£4.23.6). It is clear that both correspondences are inverse of each other and therefore they define the announced isomorphism~£4.23.1.

\bigskip
\noindent
{\bf Corollary~£4.24.} {\it  With the notation above, assume that $C$ acts transitively on $I\,.$ Then, if {\rm (Q)} holds for
 $(b_i,{}^{D_i} \!\hat K^i)$ for any $i\in I$ and  any subgroup $D_i$ of~$C_i\,,$ it holds for $(b,\hat G)\,.$\/}

\medskip
\noindent
{\bf Proof:} It follows from  Propositions~£4.21 and~£3.23 that, choosing $i\in I\,,$ we have  $\O {\rm Out}_{k^*}
(\hat G)_b\-$module isomorphisms
$$\eqalign{\R\G_k (C)\otimes_{\R\G_k (C_i)} \R_{_{\!\hat K^i\!\!}}\G_k  (\hat G^i,b_i)
&\cong \R_{_{\!\hat K\!}} \G_k  (\hat G,b)\cr 
\R\G_k (C)\otimes_{\R\G_k (C_i)}\R_{_{\K^i}}\!\G_k (\F^i,\widehat{\aut}_{(\F^i)^{^{\rm nc}}})
&\cong \R_{_{\K}}\! \G_k (\F,\widehat{\aut}_{\F^{\rm nc}})\cr}
\eqno £4.24.1.$$

\smallskip
On the other hand,  assume that we have an $\O {\rm Out}_{k^*}(\hat G^i)_{b_i}\-$module
isomorphism
$$\G_k  (\hat G^i,b_i)\cong \G_k (\F^i,\widehat{\aut}_{(\F^i)^{^{\rm nc}}})
\eqno £4.24.2;$$
then, since the restriction induces compatible $\G_k (C_i)\-$module structures on both
members of this isomorphism [10,~15.21 and~15.33], it follows from
[10,~15.23.2 and~15.37.1] that we still have an $\O {\rm Out}_{k^*}(\hat G^i)_{b_i}\-$module isomorphism 
$$\R_{_{\!\hat K_i\!\!}}\G_k  (\hat G^i,b_i)\cong \R_{_{\K^i}}\!\G_k (\F^i,\widehat{\aut}_{(\F^i)^{^{\rm nc}}})
\eqno £4.24.3.$$
Thus, from isomorphisms~£4.24.1 above, we get an $\O {\rm Out}_{k^*}(\hat G)_b\-$module isomorphism
$$\R_{_{\!\hat K\!}}\G_k  (\hat G,b)\cong \R_{_{\K}}\! \G_k (\F,\widehat{\aut}_{(\F)^{\rm nc}})
\eqno £4.24.4.$$

\smallskip
Consequently, according to our hypothesis and possibly applying Proposition~£4.18 above and [10,~15.23.2 and~15.37.1],  for any subgroup $D$ of $C$ we have an $\O {\rm Out}_{k^*}({}^D\!\hat K)_b\-$module isomorphism
$$\R_{_{\!\hat K\!}}\G_k  ({}^D\!\hat K,b)\cong \R_{_{\K}}\! \G_k ({}^D \K,\widehat{\aut}_{({}^D\K)^{\rm nc}})
\eqno £4.24.5;$$
but, since ${\rm Aut}_{k^*}(\hat G)_b$ stabilizes~$\hat K\,,$ we have evident group
homomorphisms
$$C\too {\rm Out}_{k^*}({}^{D}\!\hat K)_b\longleftarrow 
{\rm Aut}_{k^*}(\hat G)_b
\eqno £4.24.6\phantom{.}$$
and it is clear that the image of ${\rm Aut}_{k^*}(\hat G)_b$ contains and normalizes the
image of $C\,;$ hence, we still have an  $\O {\rm Out}_{k^*}(\hat G)_b\-$module 
isomorphism
$$\R_{_{\!\hat K}}\!\G_k \big({}^{D}\!\hat K,b\big)^C \cong 
\R_{_{\!\K}}\!\G_k  ({}^D \K,\widehat\aut_{({}^D \K)^{^{\rm nc}}}\big)^C
\eqno £4.24.7.$$
Then, it follows from [10,~15.23.4 and~15.38.1] that the direct sum of
isomorphisms~£4.24.6 when $D$ runs over the set of subgroups of $C$ supplies 
an $\O {\rm Out}_{k^*}(\hat G)_b\-$module isomorphism 
$\G_k (\hat G,b)\cong \G_k (\F,\widehat\aut_{\F^{^{\rm nc}}})\,.$ We are done.
\vfill
\eject

\bigskip
\bigskip
\noindent
{\bf References}
\bigskip
\noindent
[1]\phantom{.} Michel Brou\'e and Llu\'\i s Puig, {\it Characters and Local
Structure in $G\-$alge-bras,\/} Journal of Algebra, 63(1980), 306-317.
\smallskip\noindent
[2]\phantom{.} Yun Fan and Lluis Puig, {\it On blocks with nilpotent
coefficient extensions\/}, Algebras and Representation Theory, 1(1998),
27-73 and Publisher revised form, 2(1999), 209.
\smallskip\noindent
[3]\phantom{.} Daniel Gorenstein, {\it ``Finite groups''\/}, Harper's Series,
1968, Harper and Row.
\smallskip\noindent
[4]\phantom{.} James Green, {\it Some remarks on defect groups\/}, Math. Zeit.,
107(1968), 133-150.
\smallskip\noindent
[5]\phantom{.} Llu\'\i s Puig, {\it Pointed groups and  construction of
characters}, Math. Zeit. 176(1981), 265-292. 
\smallskip\noindent
[6]\phantom{.} Llu\'\i s Puig, {\it Local fusions in block source algebras\/},
Journal of Algebra, 104(1986), 358-369. 
\smallskip\noindent
[7]\phantom{.} Llu\'\i s Puig, {\it Nilpotent blocks and their source
algebras}, Inventiones math., 93(1988), 77-116.
\smallskip\noindent
[8]\phantom{.} Llu\'\i s Puig, {\it Pointed groups and  construction of
modules}, Journal of Algebra, 116(1988), 7-129.
\smallskip\noindent
[9]\phantom{.} Llu\'\i s Puig, {\it ``Blocks of Finite Groups''\/},
Springer Monographs in Mathematics, 2002, Springer-Verlag, Berlin, Barcelona.
\smallskip\noindent
[10]\phantom{.} Llu\'\i s Puig, {\it ``Frobenius categories versus Brauer blocks''\/}, Progress in Math., 
274(2009), Birkh\"auser, Basel.
\smallskip\noindent
[11]\phantom{.} Llu\'\i s Puig, {\it Block Source Algebras  in p-Solvable
Groups},  Michigan Math. J. 58(2009), 323-328
\smallskip\noindent
[12]\phantom{.} Llu\'\i s Puig, {\it Nilpotent extension of blocks},  submitted to Math. Z.

\bigskip
\bigskip
\noindent
{\bf Abstract}
\medskip
\noindent
We show that the refinement of Alperin's Conjecture proposed in [10, Ch.~16] can be proved by checking that this
refinement holds on any central $k^*\-$exten-sion of a finite group $H$ containing a normal simple group $S$ with trivial
centralizer in $H$ and $p'\-$cyclic quotient $H/S\,.$ This paper improves our result in [10,~Theorem~16.45] and repairs
some bad arguments there.

\end